\newif\ifAIAA\AIAAfalse
\newcommand{\thetitle}{Optimized Runge--Kutta (LDDRK) timestepping schemes for non-constant-amplitude oscillations}
\newcommand{\prevpub}{A preliminary version of some parts of this paper was presented as part of AIAA Paper 2019-2414 at the 25th AIAA/CEAS Aeroacoustics Conference in Delft, The Netherlands~\citep{petronilia+brambley-2019-aiaa}.}

\ifAIAA
  \PassOptionsToPackage{unicode}{hyperref}
  \documentclass[conf]{new-aiaa}
  \usepackage{fancyhdr}
\else
  \PassOptionsToPackage{unicode}{hyperref}
  \documentclass[3p,number,sort&compress]{elsarticle}
  \journal{Journal of Computational Physics}
  \biboptions{longnamesfirst}
  
  \makeatletter
    \providecommand{\doi}[1]{%
      \begingroup
        \let\bibinfo\@secondoftwo
       \urlstyle{rm}%
        \href{http://dx.doi.org/#1}{%
          doi:\discretionary{}{}{}%
          \nolinkurl{#1}%
        }%
      \endgroup
    }
  \makeatother
\fi

\usepackage[utf8]{inputenc}
\usepackage{graphicx}
\usepackage{amsmath}
\usepackage{psfrag}
\usepackage{color}
\usepackage{subcaption}
\usepackage{xcolor}
\newcommand{\eps}{\varepsilon} % \eps: Epsilon
\newcommand{\vect}[1]{\boldsymbol{#1}} % \vect: Vectors
\newcommand{\Real}{\mathrm{Re}} % \Real: Real part of
\newcommand{\I}{\mathrm{i}} % \I: $\sqrt{-1}$
\newcommand{\e}{\mathrm{e}} % \e: $\exp{1}$
\newcommand{\intd}{\mathrm{d}} % \intd: Roman d for integrals, eg. \int\ldots\,\intd\theta
\newcommand{\odt}{\omega\Delta t}
\newcommand{\obdt}{\bar{\omega} \Delta t}
\definecolor{yellow}{HTML}{F4E824}
\definecolor{orange}{HTML}{FEC338}
\definecolor{green}{HTML}{18BFB5}

\def\clap#1{\hbox to 0pt{\hss#1\hss}}

\def\mathclap{\mathpalette\mathclapinternal}

\def\mathclapinternal#1#2{%
\clap{$\mathsurround=0pt#1{#2}$}}

\ifAIAA
  \title{\thetitle}

  \author{Alda\"\i r Petronilia\footnote{Undergraduate, University of Warwick} and Edward James Brambley\footnote{Associate Professor, University of Warwick, AIAA senior member.}}
  \affil{University of Warwick, Coventry CV4 7AL, United Kingdom}

  \fancypagestyle{firststyle}
  {
    \fancyhf{}
    
    \lhead{25th AIAA/CEAS Aeroacoustics Conference, Delft, The Netherlands, May 2019.}
    \rhead{AIAA Paper 2019--2414}
    \cfoot{\thepage}
  }

\fi

\begin{document}

\ifAIAA
  \maketitle
  \thispagestyle{firststyle}

\else
  % For JCP

  \begin{frontmatter}

  \title{\thetitle\tnoteref{t1}}
  \tnotetext[t1]{\prevpub}

  \author[maths]{Alda\"\i r~Petronilia}
  \author[maths,wmg]{Edward~James~Brambley\corref{cor1}}
  \ead{E.J.Brambley@warwick.ac.uk}

  \cortext[cor1]{Corresponding author.  Tel.: +44 2476 574832.}

  \address[maths]{Mathematics Institute, University of Warwick, Coventry, CV4 7AL, United Kingdom}
  \address[wmg]{WMG, University of Warwick, Coventry, CV4 7AL, United Kingdom}
\fi

%%%%%%%%%%%%
% Document %
%%%%%%%%%%%%

\begin{abstract}
  Finite differences and Runge--Kutta time stepping schemes used in Computational AeroAcoustics simulations are often optimized for low dispersion and dissipation (e.g.\ DRP or LDDRK schemes) when applied to linear problems in order to accurately simulate waves with the least computational cost.  Here, the performance of optimized Runge--Kutta time stepping schemes for linear time-invariant problems with non-constant-amplitude oscillations is considered.  This is in part motivated by the recent suggestion that optimized spatial derivatives perform poorly for growing and decaying waves, as their optimization implicitly assumes real wavenumbers.  To our knowledge, this is the first time the time-stepping of non-constant-amplitude oscillations has been considered.  It is found that current optimized Runge--Kutta schemes perform poorly in comparison with their maximal order equivalents for non-constant-amplitude oscillations.  Moreover, significantly more accurate results can be achieved for the same computation cost by replacing a two-step scheme such as LDDRK56 with a single step higher-order scheme with a longer time step.
  Attempts are made at finding optimized schemes that perform well for non-constant-amplitude oscillations, and three such examples are provided.  However, the traditional maximal order Runge--Kutta time stepping schemes are still found to be preferable for general problems with broadband excitation.
  These theoretical predictions are illustrated using a realistic 1D wave-propagation example.
\end{abstract}
  
\ifAIAA\else\end{frontmatter}\fi

\section{Introduction}
Computational AeroAcoustics (CAA) simulations are an important tool in investigating aircraft noise in realistic geometries and flows.  Unlike Computational Fluid Dynamics (CFD) simulations, CAA simulations are designed to accurately propagate small amplitude oscillations over the entire computational domain, and are therefore poised on a knife edge between being overly dissipative on the one hand and being unstable on the other.  Finite difference schemes to calculate spatial derivatives, and Runge--Kutta and Adams--Bashforth schemes to step forwards in time, have all been optimized to attempt to accurately propagate acoustic perturbations with few points per wavelength and steps per period respectively, and such schemes are commonly used in modern CAA simulations.  Examples of optimized spatial derivatives include: the by-now classic 7-point 4th order explicit DRP schemes~\citep{tam+webb-1993,tam+shen-1993}; optimized implicit/compact schemes of up to 6th order~\citep{kim+lee-1996,kim-2007}; prefactored implicit MacCormack schemes~\citep{hixon+turkel-2000}; trigonometrically optimized schemes~\citep{tang+baeder-1999}; 2nd and 4th order 9, 11, and 13 point schemes~\citep{bogey+bailly-2004}; and asymmetric optimized schemes for use near boundaries~\citep{berland+bogey+marsden+bailly-2007,turner+haeri+kim-2016}.  Examples of optimized timestepping schemes include: an optimized Adams--Bashforth scheme~\citep{tam+webb-1993}; Low Dispersion and Dissipation Runge--Kutta (LDDRK) 5-stage, 6-stage and alternating 4/6- and 5/6-stage schemes~\citep{hu+hussaini+manthey-1996}; and optimized 5- and 6-stage Runge--Kutta schemes~\citep{bogey+bailly-2004}.  \Citet{rona+etal-2017} even considered jointly optimizing spatial derivatives and time stepping schemes to give the best wave propagation properties when combined, although this analysis is dependent on the dispersion relation of the system being simulated, while the previously mentioned schemes are applicable to general dispersion relations.  Such optimized schemes are also starting to find favour outside of computational aeroacoustics, for example in Large Eddy Simulations~\citep[e.g.][]{sciacovelli+cinnella+grasso-2017,cohen+gloerfelt-2018,capuano+bogey+spelt-2018}

Most of the optimized spatial and time-stepping schemes investigate the action of the scheme in the wavenumber/frequency domain.  For example, the solution to the time-stepping problem $\intd\vect{U}\!/\intd t = \vect{F}(\vect{U}, t)$ for the $p$-stage low-storage Runge--Kutta scheme considered by~\citet*{hu+hussaini+manthey-1996} is
\begin{align}
  \vect{U}(t+\Delta t) &= \vect{U}(t) + \beta_p\vect{K_p},&
  &\text{where}\qquad \vect{K_{j+1}} = \Delta t\vect{F}\big(\vect{U}(t) + \beta_j\vect{K}_j,\, t+\beta_j\Delta t\big),
\end{align}
with $\beta_0 = 0$.  Assuming $F(\vect{U},t)$ to be linear in $\vect{U}$ and time invariant, and transforming to the frequency domain (or equivalently considering $F(\vect{U}, t) = -\I\omega\vect{U}$), this scheme results in
\begin{align}
  \vect{U}(t+\Delta t) &= r(\odt)\vect{U}(t),&
  &\text{where}\qquad r(\odt) = 1 + \sum_{j=1}^p c_j(-\I\odt)^j &
  &\text{and}\qquad \beta_{p-j} = c_{j+1}/c_j
\label{equ:r}
\end{align}
with coefficient $c_1=\beta_p$.  Since the exact solution would have $\vect{U}(t+\Delta t) = \vect{U}(t)r_e(\odt)$ with $r_e(\odt) = \exp\{-\I\odt\}$, we may define the effective numerical angular frequency $\bar{\omega}$ by $r(\odt) = \exp\{-\I\obdt\}$.  One could choose the coefficients $c_j = 1/j!$ so that $|r(\odt)-r_e(\odt)| = O\big((\Delta t)^{p+1}\big)$ in the limit $\Delta t\to 0$; this is referred to as a $p$th-order accurate scheme, and since it is the best that can be achieved with a $p$-stage Runge--Kutta scheme, we refer to such schemes here as maximal order.  In contrast, one could instead vary the coefficients $c_j$ in order to minimize an error of the form
\begin{align}
  e &= \int_0^{\pi\eta} \left|r(\odt) - r_e(\odt)\right|^2\,\intd(\odt),&
  &\text{or}&
  E &= \int_0^{\pi\eta} \left|\bar\odt - \odt\right|^2\,\intd(\odt),
\label{equ:metric}
\end{align}
subject to constraints of a minimum order of accuracy (typically 2nd or 4th order accuracy as $\Delta t\to 0$) and stability (meaning $|r|\leq 1$ for $0\leq \odt < \pi\eta_s$ for some given stability threshold $\eta_s$).  Optimization of $e$ was performed by~\citet{hu+hussaini+manthey-1996}, while \citet{tam+webb-1993} optimized the equivalent of $E$ for an Adams--Bashforth scheme\footnote{Note that~\citet{tam+webb-1993} used the opposite notation to that used here, in that they used $(\bar\omega, \omega)$ where here we use $(\omega, \bar\omega)$.}.  A similar optimization method may be performed for spatial derivative schemes in terms of the spatial wavenumber $k\Delta x$ instead of the frequency $\odt$~\citep[see, for example,][]{tam+webb-1993}.

Recently~\citep{brambley-2016-jcp}, it was suggested for spatial derivatives that optimization of a metric such as~\eqref{equ:metric} which assumes real $k\Delta x$ results in a scheme which performs well for constant amplitude waves corresponding to real $k$, but which performs poorly for waves of non-constant amplitude corresponding to complex $k$.  Unfortunately, non-constant amplitude waves are rather common in aeroacoustics, especially in the vicinity of acoustic linings, for high-order spinning modes which decay rapidly away from duct walls, close to near-singularities such as sharp trailing edges or strongly localized sources, and for instabilities;  non-constant-amplitude oscillations are also common in other branches of physics.  Attempts at reoptimizing spatial derivatives to perform well for both non-constant and constant amplitude waves~\citep{brambley+markeviciute-2017-aiaa} concluded that, with sufficient a priori knowledge of expected wavenumbers, optimized derivative schemes could be constructed to perform well, but that in general, and certainly for broadband excitation, maximal order schemes were more likely to be more accurate.

The purpose of this paper is to investigate the comparable situation for timestepping schemes.  In particular, we consider how well common schemes (such as the LDDRK56 scheme~\citep{hu+hussaini+manthey-1996}) which are optimized using metrics such as~\eqref{equ:metric} perform when $\odt$ is not real, corresponding to waves which decay or grow exponentially in time.  In section~\ref{section:TheoreticalComparison} we consider the theoretical performance of Runge--Kutta schemes in the frequency domain, and compare the accuracy and stability of a number of existing schemes for complex frequencies.  Informed by this, in section~\ref{section:optimization} we reconsider optimizing Runge--Kutta schemes for oscillations of non-constant amplitude.  All such schemes are compared in practice using a simple 1D wave propagation test case in section~\ref{section:1d}.  Finally, conclusion about the performance of Runge--Kutta schemes for non-constant-amplitude oscillations are given in section~\ref{section:conclusion}, including recommendations for which schemes to use in practical applications and potential avenues for future research.

\section{Theoretical comparison of Runge--Kutta time stepping schemes}
\label{section:TheoreticalComparison}

Any Runge--Kutta scheme solves the linear differential equation $\intd\vect{U}\!/\intd t = -\I\omega\vect{U}$ to give $\vect{U}(t+\Delta t) = \vect{U}(t)r(\odt)$, where $r(\odt)$ is given by~\eqref{equ:r} with some constants $c_j$.  We define the relative phase error of such a scheme by
\begin{align}
  \eps_p &= \left|\frac{\bar\omega - \omega}{\omega}\right| = \left|\frac{\bar\omega}{\omega} - 1 \right|,&
  &\text{where}\qquad \bar\odt = \I\log\big(r(\odt)\big),
\end{align}
and the branch of the logarithm is chosen so that $|\bar\odt-\odt|$ is minimized.
We may also define the relative amplification factor error
\begin{equation}
\eps_r = \left|\frac{r(\odt)-r_e(\odt)}{r_e(\odt)}\right| = \big|r(\odt)\exp\{\I\odt\}-1\big|.
\end{equation}
Assuming these errors are small, the relative global error for a simulation up to time $t=T$ compared with the exact solution $\vect{U_e}(t)$ is then given by
\begin{equation}
  \frac{\|\vect{U}(t) - \vect{U_e}(t)\|}{\|\vect{U_e}(t)\|}
  = \left|\frac{r(\odt)^{T/\Delta t}-r_e(\odt)^{T/\Delta t}}{r_e(\odt)^{T/\Delta t}}\right|
  \approx \frac{T}{\Delta t}\eps_r
  \approx T \left |\omega\right |\eps_p.
\label{equ:accuracy}
\end{equation}
From this, we may interpret the relative amplification factor error $\eps_r$ as being the time-stepping error per time step, since there are $T/\Delta t$ timesteps in the whole simulation.  We may similarly interpret the relative phase error $\eps_p$ as being the time-stepping error per oscillation period, up to a factor of $2\pi$, at least for real $\omega$, since there are $T|\omega|/2\pi$ oscillation periods in the whole simulation.  Since the number of oscillation periods in the whole simulation is independent of the numerical timestep $\Delta t$ used, we will adopt the relative phase error $\eps_p$ as our measure of accuracy in what follows.

\subsection{Phase errors for complex frequencies}

\begin{figure}
  \begin{subfigure}[b]{0.32\textwidth}
    \includegraphics[width=\textwidth]{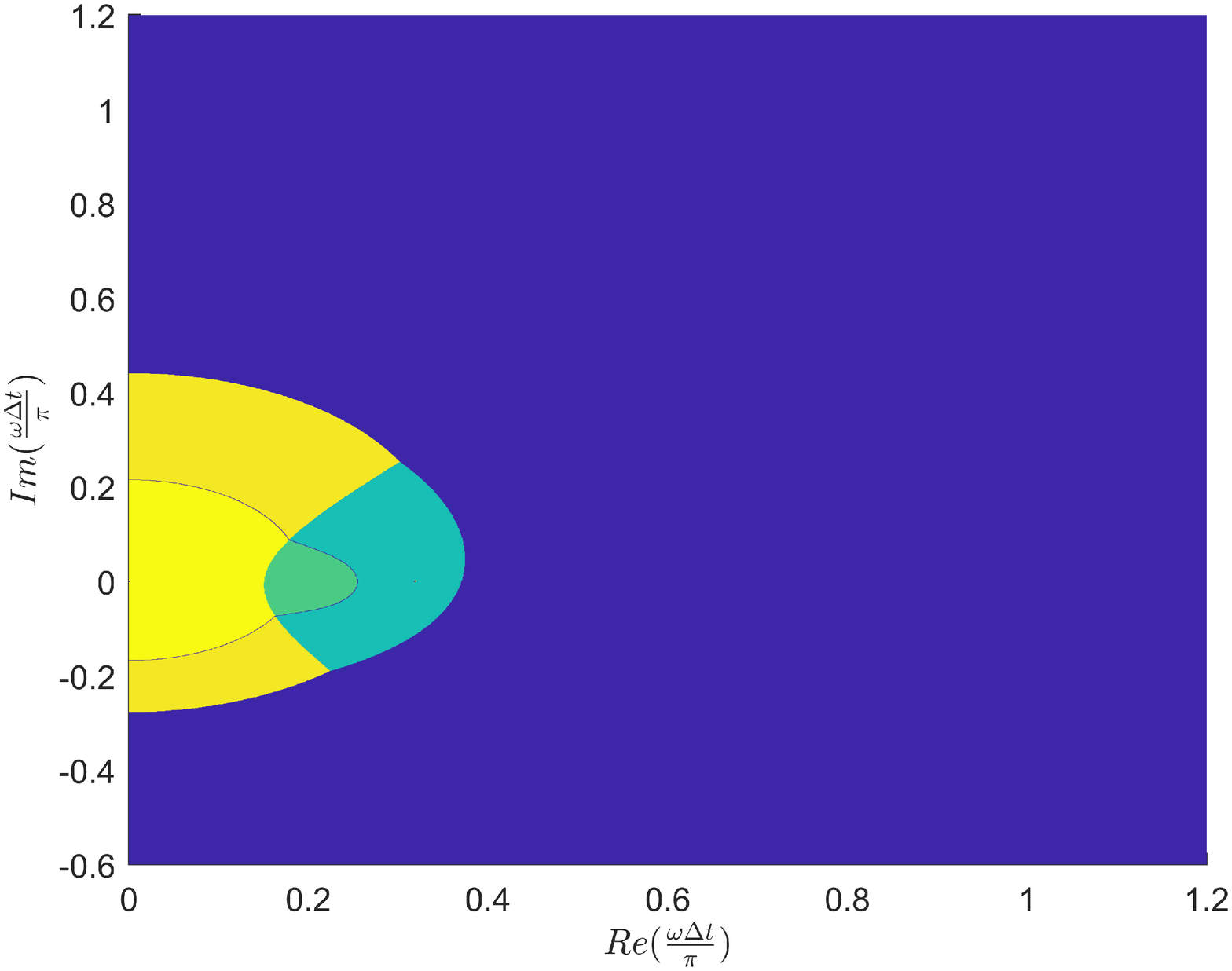}
    \caption{\textcolor{yellow}{4th Order} vs \textcolor{green}{LDDRK4}}
    \label{fig:LDDRK4VSRK4OnePercentError}
  \end{subfigure}
  \begin{subfigure}[b]{0.32\textwidth}
    \includegraphics[width=\textwidth]{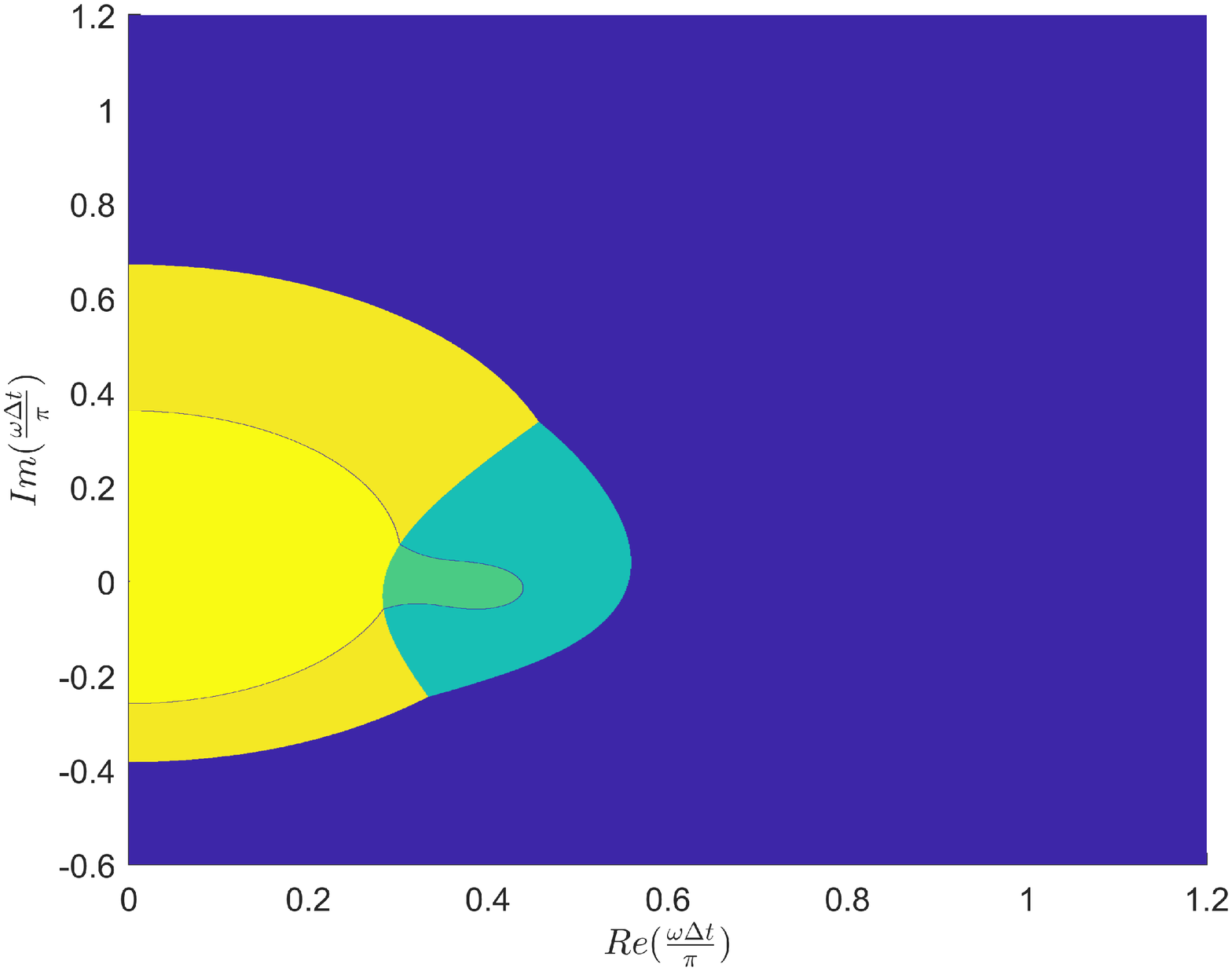}
    \caption{\textcolor{yellow}{5th Order} vs \textcolor{green}{LDDRK5}}
    \label{fig:LDDRK5VSRK6OnePercentError}
  \end{subfigure}
  \begin{subfigure}[b]{0.32\textwidth}
    \includegraphics[width=\textwidth]{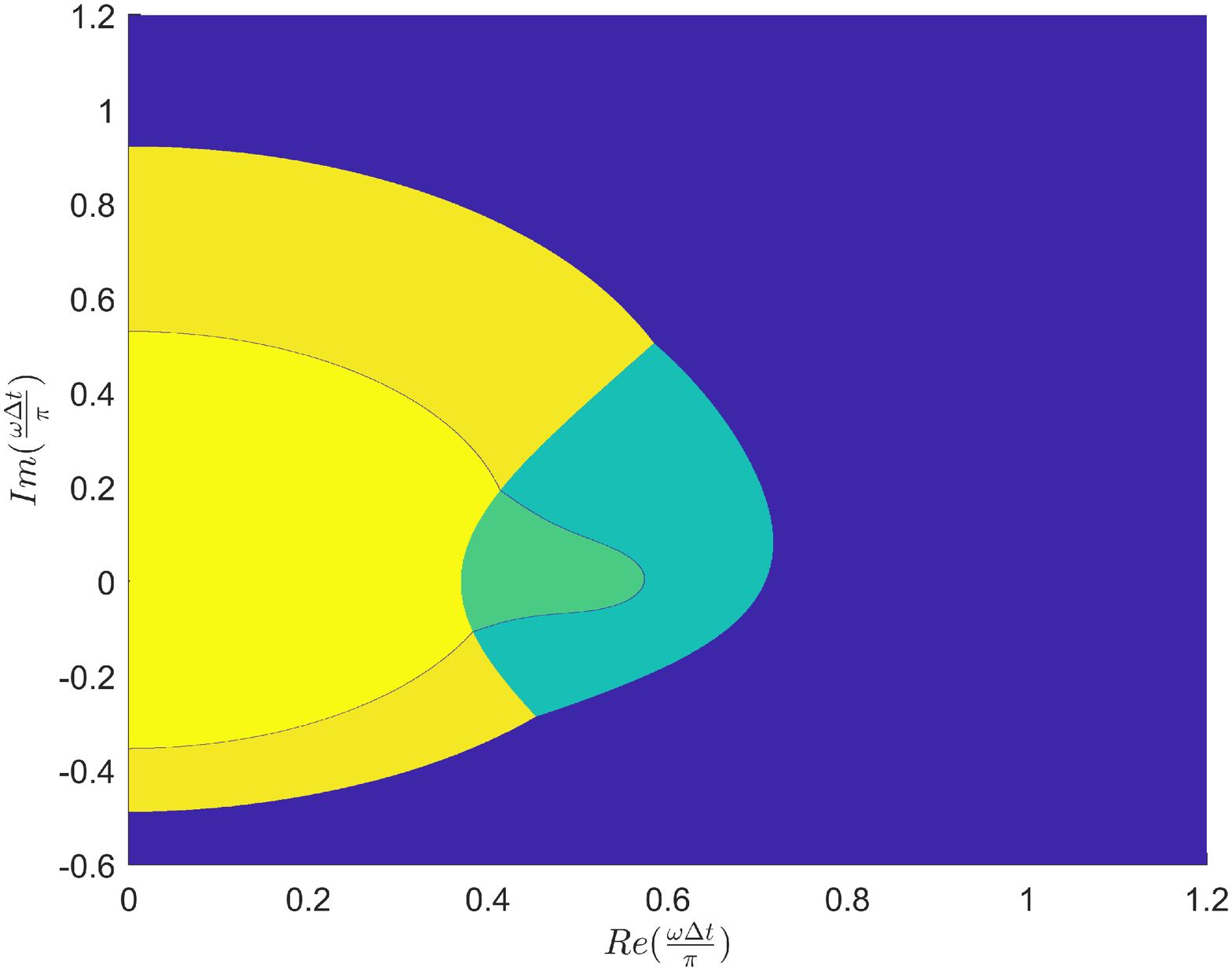}
    \caption{\textcolor{yellow}{6th Order} vs \textcolor{green}{LDDRK6}}
    \label{fig:LDDRK6VSRK6OnePercentError}
  \end{subfigure}
  \caption{Comparison in the complex $\odt$ plane of where the maximal-order (yellow) or LDDRK (green) schemes are more accurate for the 4, 5 and 6 stage schemes.  Darker regions indicate where neither scheme is $0.1\%$ accurate ($\eps_p > 10^{-3}$), while violet regions indicate where neither scheme is $1\%$ accurate ($\eps_p > 10^{-2}$).}
  \label{fig:LDDRKvsRkOnePercentPhaseError}
\end{figure}

We first consider the single-step Low Dispersion and Dissipation Runge--Kutta (LDDRK) schemes of~\citet{hu+hussaini+manthey-1996}.
\begin{figure}
  \begin{subfigure}[b]{0.49\textwidth}
    \includegraphics[width=\textwidth]{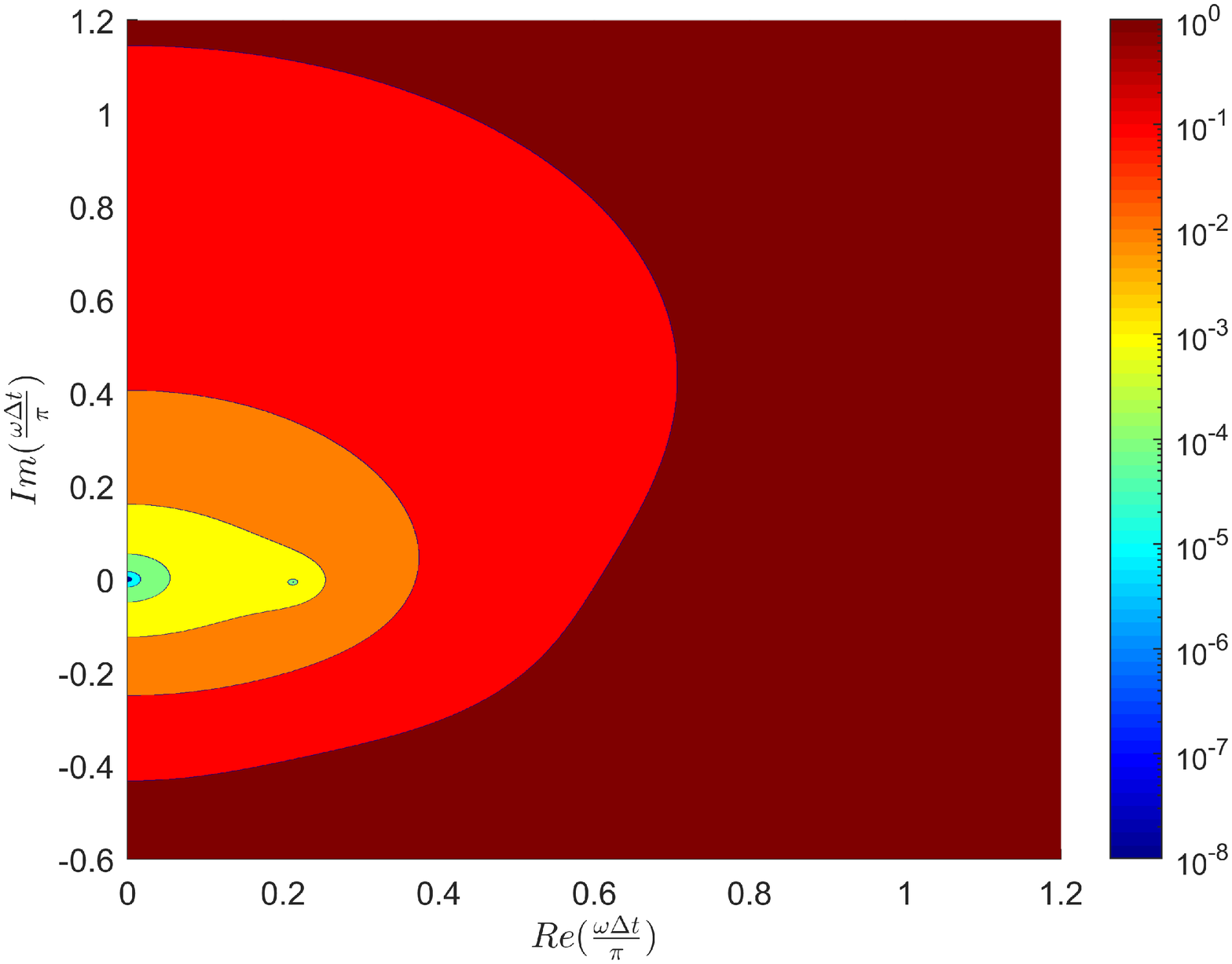}
    \caption{2nd Order LDDRK4}
    \label{fig:LDDRK4}
  \end{subfigure}
  \hfill
  \begin{subfigure}[b]{0.49\textwidth}
    \includegraphics[width=\textwidth]{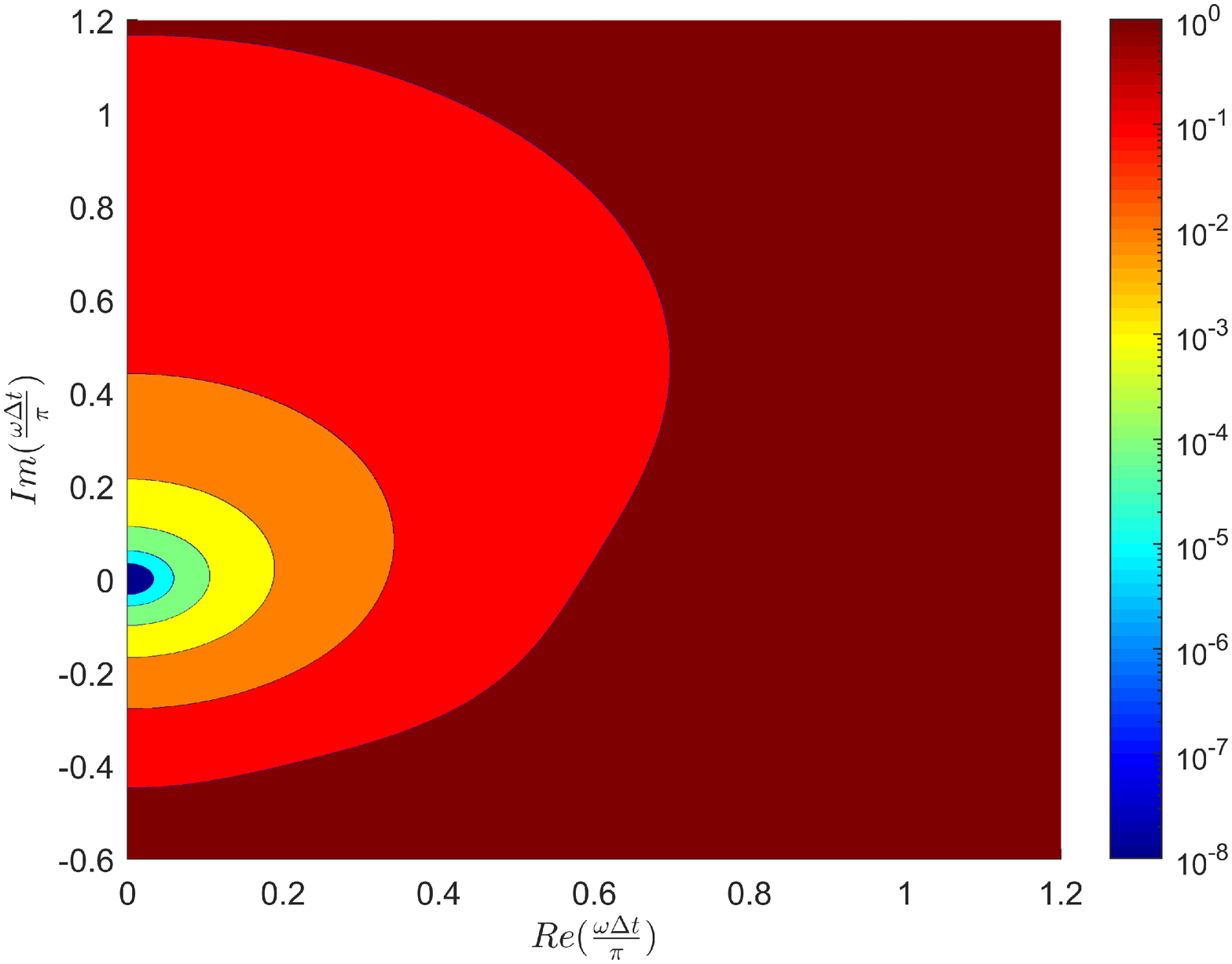}
    \caption{4th Order RK4}
    \label{fig:RK4}
  \end{subfigure}
  
  \begin{subfigure}[b]{0.49\textwidth}
    \includegraphics[width=\textwidth]{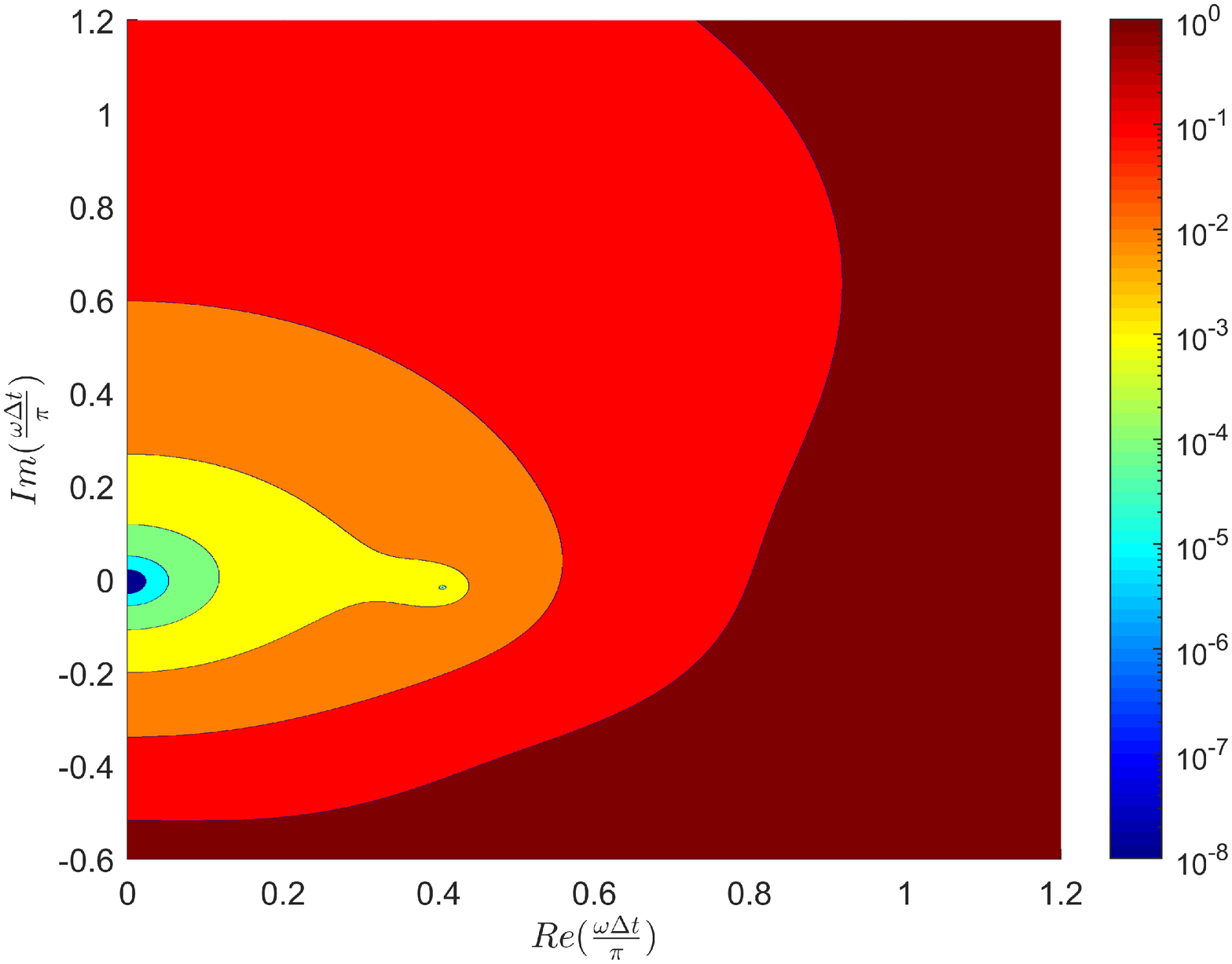}
    \caption{2nd Order LDDRK5}
    \label{fig:LDDRK5}
  \end{subfigure}
  \hfill
  \begin{subfigure}[b]{0.49\textwidth}
    \includegraphics[width=\textwidth]{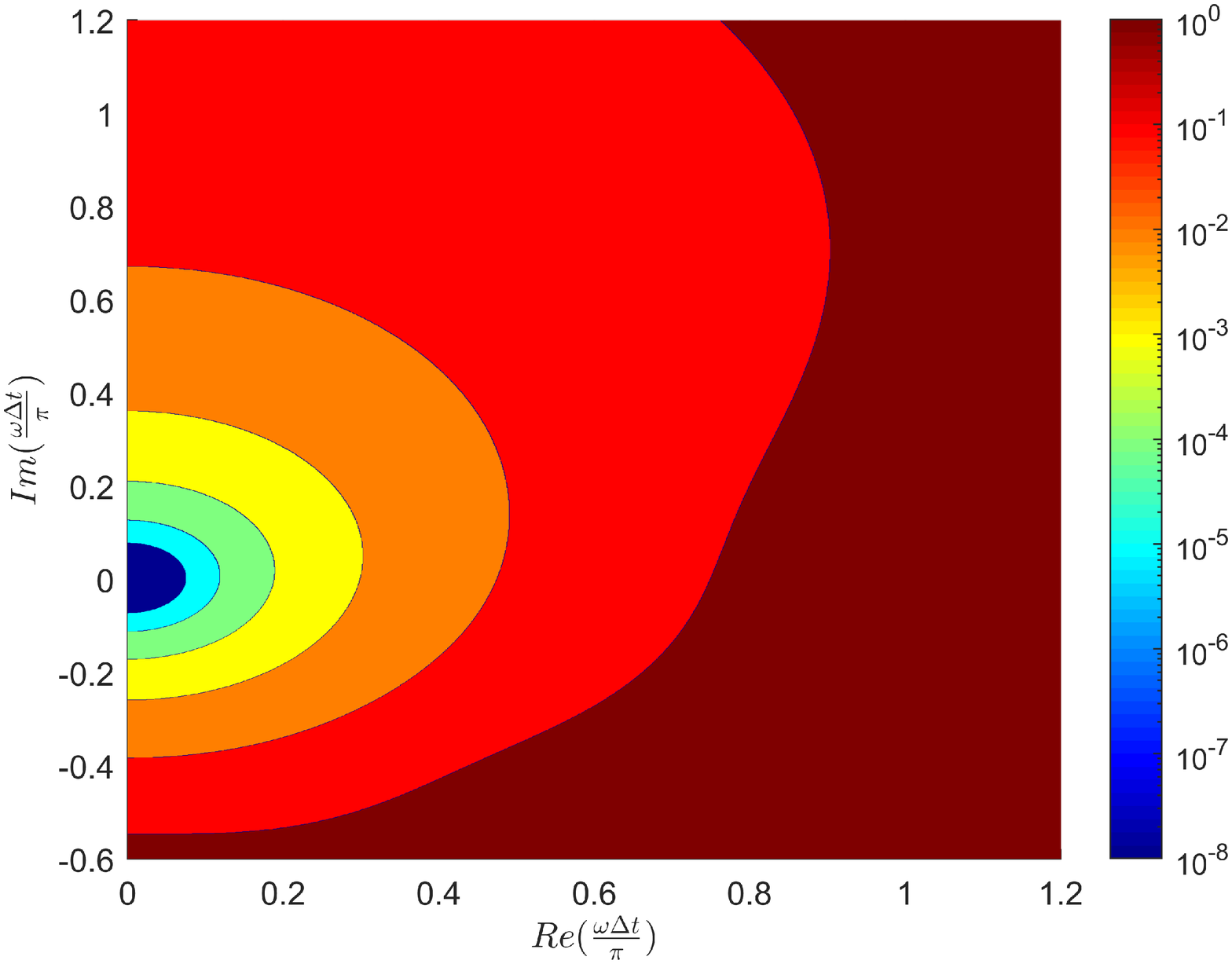}
    \caption{5th order RK5}
    \label{fig:RK5}
  \end{subfigure}
  
  \begin{subfigure}[b]{0.49\textwidth}
    \includegraphics[width=\textwidth]{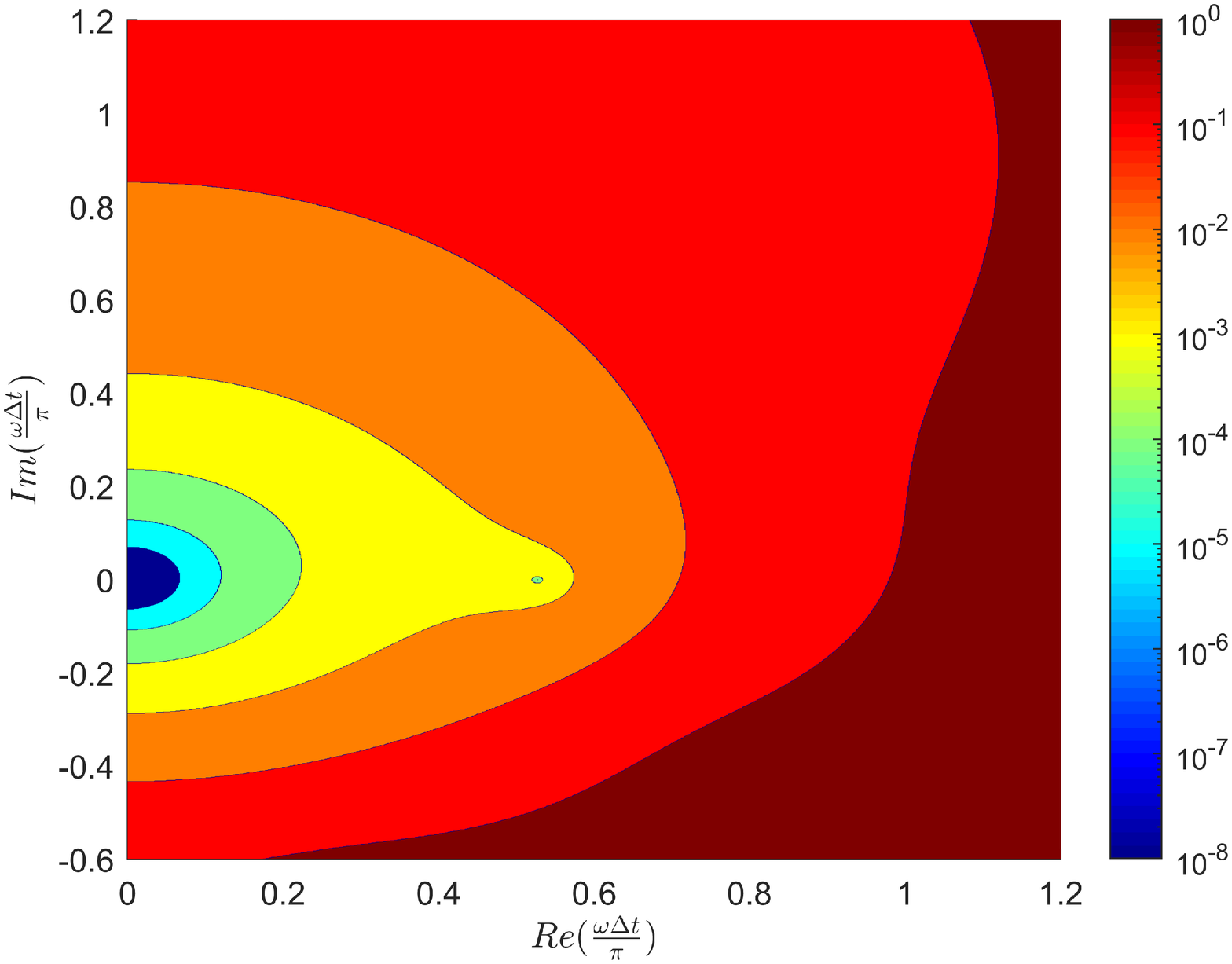}
    \caption{4th order LDDRK6}
    \label{fig:LDDRK6}
  \end{subfigure}
  \hfill
  \begin{subfigure}[b]{0.49\textwidth}
    \includegraphics[width=\textwidth]{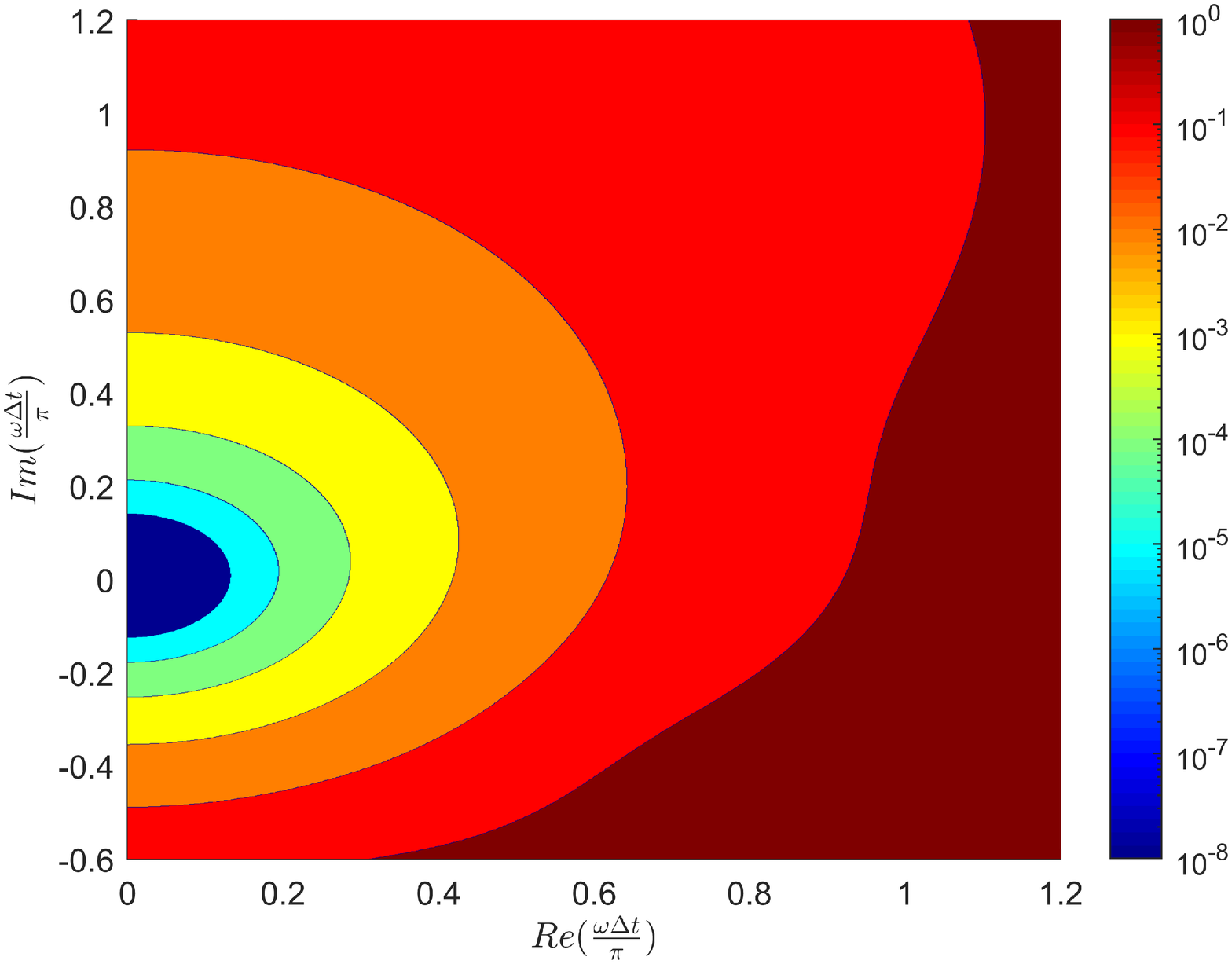}
    \caption{6th order RK6}
    \label{fig:RK6}
  \end{subfigure}
  \caption{Plots of phase error $\eps_p$ in the complex $\odt$ plane. (a) and (b) are 4-stage, (c) and (d) are 5-stage, and (e) and (f) are 6-stage. (a) and (c) are optimised 2nd order schemes~\citep{hu+hussaini+manthey-1996}, (e) is an optimized 4th order scheme~\citep{hu+hussaini+manthey-1996}, while (b), (d), and (f) are maximal order schemes.}
  \label{fig:LDDRKVSRK}
\end{figure}
Figure~\ref{fig:LDDRKVSRK} compares the phase error $\eps_p$ of these optimized LDDRK schemes with their maximal order Runge--Kutta equivalents with the same number of stages, and hence with the same computational expense.  As expected, the optimized schemes perform better along the real $\odt$ axis, at the expense of their behaviour near the origin and for non-real $\odt$.  In order to better compare the schemes, it is helpful to plot which is more accurate for any given value of $\odt$, shown in figure~\ref{fig:LDDRKvsRkOnePercentPhaseError}.
The green regions show where the LDDRK optimized schemes outperform the maximal order schemes, but operating within such a limited region would need significant a priori knowledge of the behaviour expected, and remaining in this region would be nearly impossible for broadband simulations containing a wide range of frequencies.

\begin{figure}
  \begin{subfigure}[b]{0.42\textwidth}
    \includegraphics[width=\textwidth]{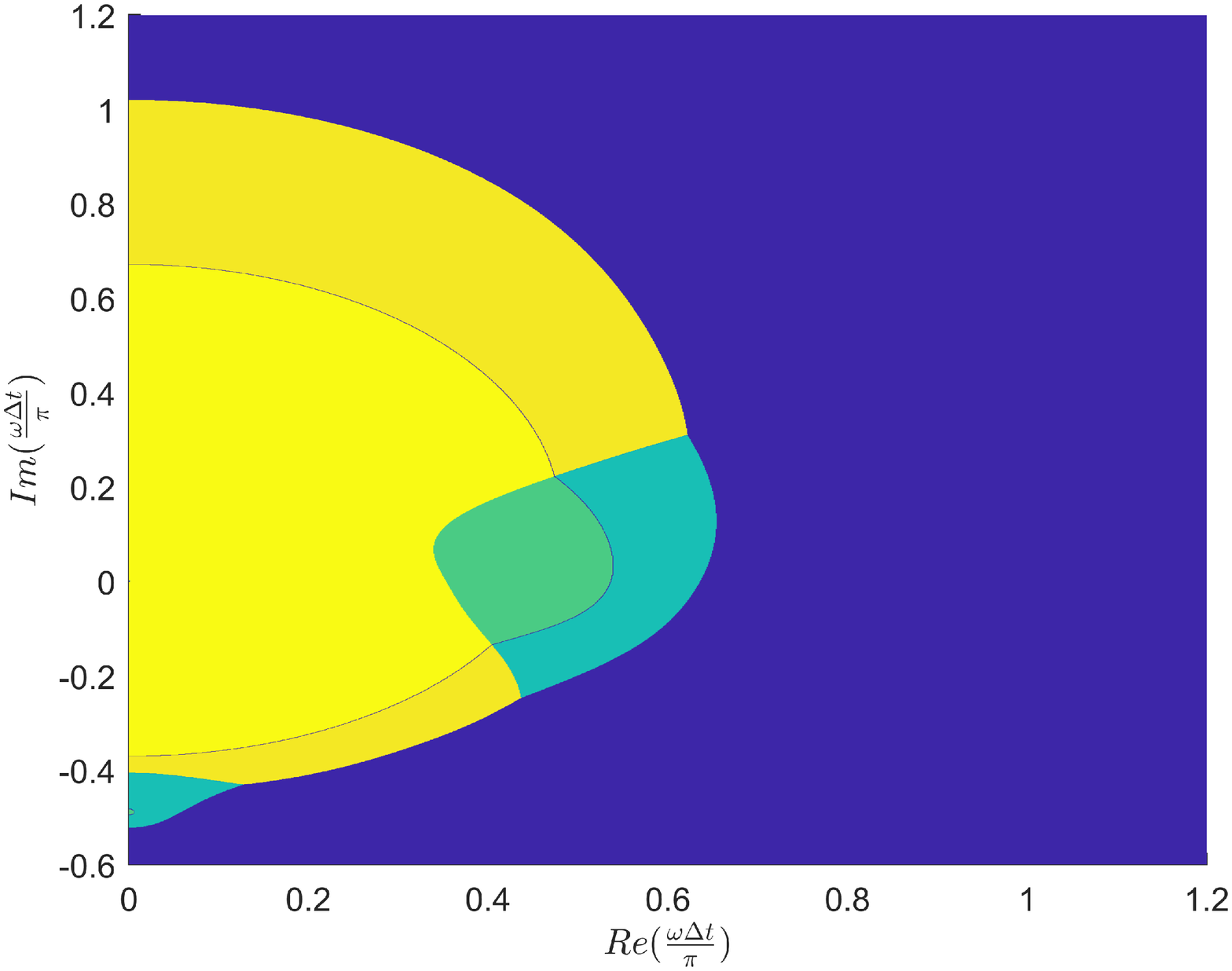}
    \caption{\textcolor{yellow}{RK10($2\Delta t$)} vs \textcolor{green}{LDDRK46($\Delta t$)}}
    \label{fig:LDDRK46VSRK10OnePercentError}
  \end{subfigure}
  \hfill
  \begin{subfigure}[b]{0.42\textwidth}
    \includegraphics[width=\textwidth]{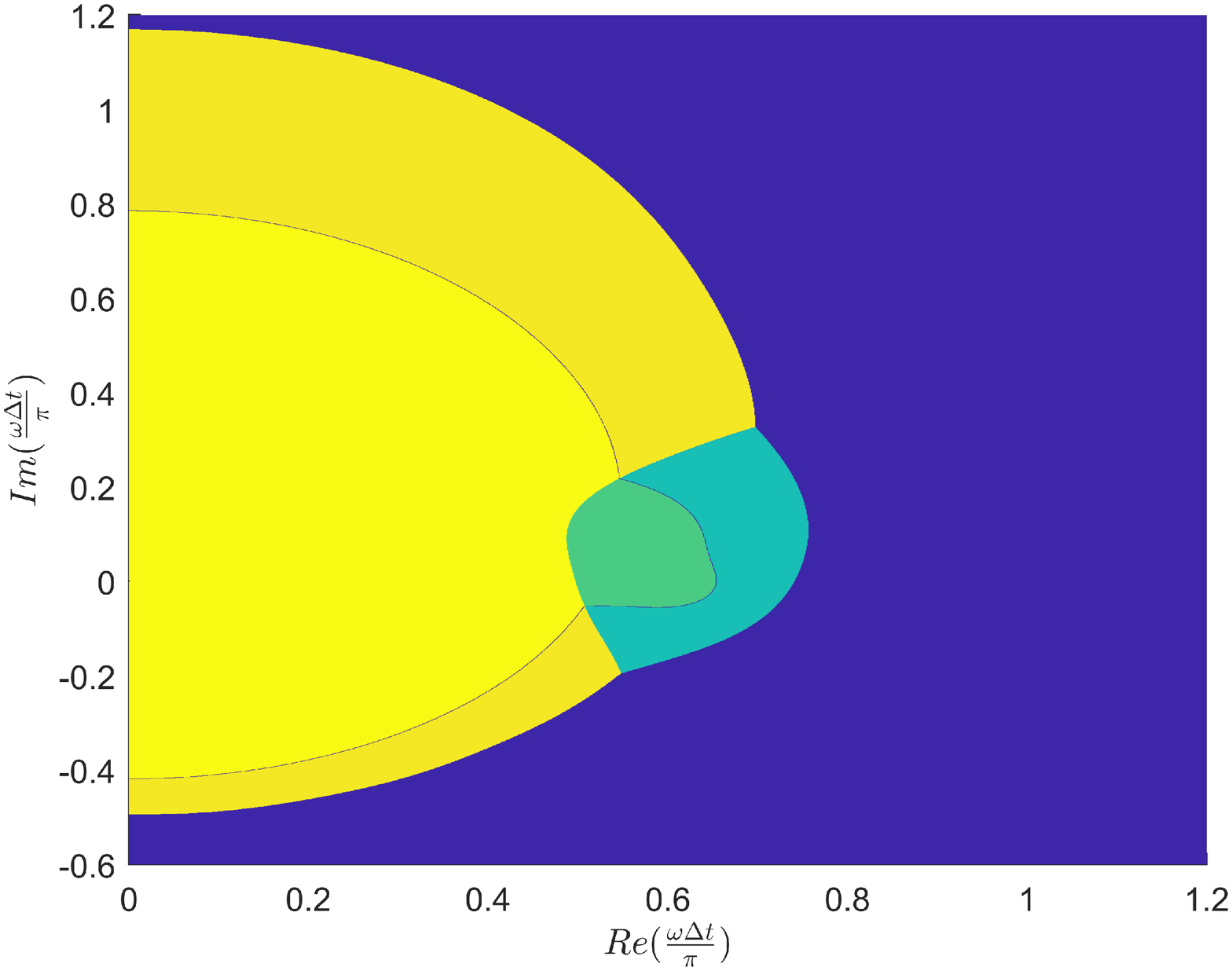}
    \caption{\textcolor{yellow}{RK11($2\Delta t$)} vs \textcolor{green}{LDDRK56($\Delta t$)}}
    \label{fig:LDDRK56VSRK11OnePercentError}
  \end{subfigure}
  \caption{Comparison in the complex $\odt$ plane of where the maximal-order (yellow) or LDDRK (green) schemes are more accurate.  Darker regions indicate where neither scheme is $0.1\%$ accurate ($\eps_p > 10^{-3}$), while violet regions indicate where neither scheme is $1\%$ accuracy ($\eps_p > 10^{-2}$).}
  \label{fig:twoStepLDDRKvstwoTimeStepRkOnePercentPhaseError}
\end{figure}

\Citet{hu+hussaini+manthey-1996} further proposed two two-step alternating schemes.  For a $p_1$-$p_2$ scheme, the first timestep uses a $p_1$-stage Runge-Kutta scheme and the second timestep uses a $p_2$-stage Runge-Kutta scheme.  This results in $\vect{U}(t + 2\Delta t) = \vect{U}(t)r_1(\odt)r_2(\odt)$, where the amplification factors of the first and second steps are
\begin{align}
    r_1(\odt) &= 1 + \sum_{j = 1}^{p_1} a_j(-\I\odt)^j = \e^{-\I\bar\omega_1\Delta t}, &
    r_2(\odt) &= 1 + \sum_{j = 1}^{p_2} b_j(-\I\odt)^j = \e^{-\I\bar\omega_2\Delta t}.
\end{align}
The target optimization in this case is given by
\begin{equation}
    \int_0^{\pi\eta} \left| r_1(\odt)r_2(\odt)  - \e^{-2\I\odt} \right|^2 \intd(\odt),
\end{equation}
while the effective numerical frequency is $\bar\omega = (\bar\omega_1 + \bar\omega_2)/2$, which may be used to calculate the phase error $\eps_p$.
Since a $p_1$-$p_2$ scheme evaluates $\vect{F}(\vect{U})$ $p_1+p_2$ times to step the time forward by $2\Delta t$, such schemes have the same computational cost as a single step of a ($p_1+p_2$)-stage Runge--Kutta scheme with a time step of $2\Delta t$.  \Citet{hu+hussaini+manthey-1996} proposed a 4-6 and a 5-6 alternating scheme.
\begin{figure}
  \begin{subfigure}[b]{0.49\textwidth}
    \includegraphics[width=\textwidth]{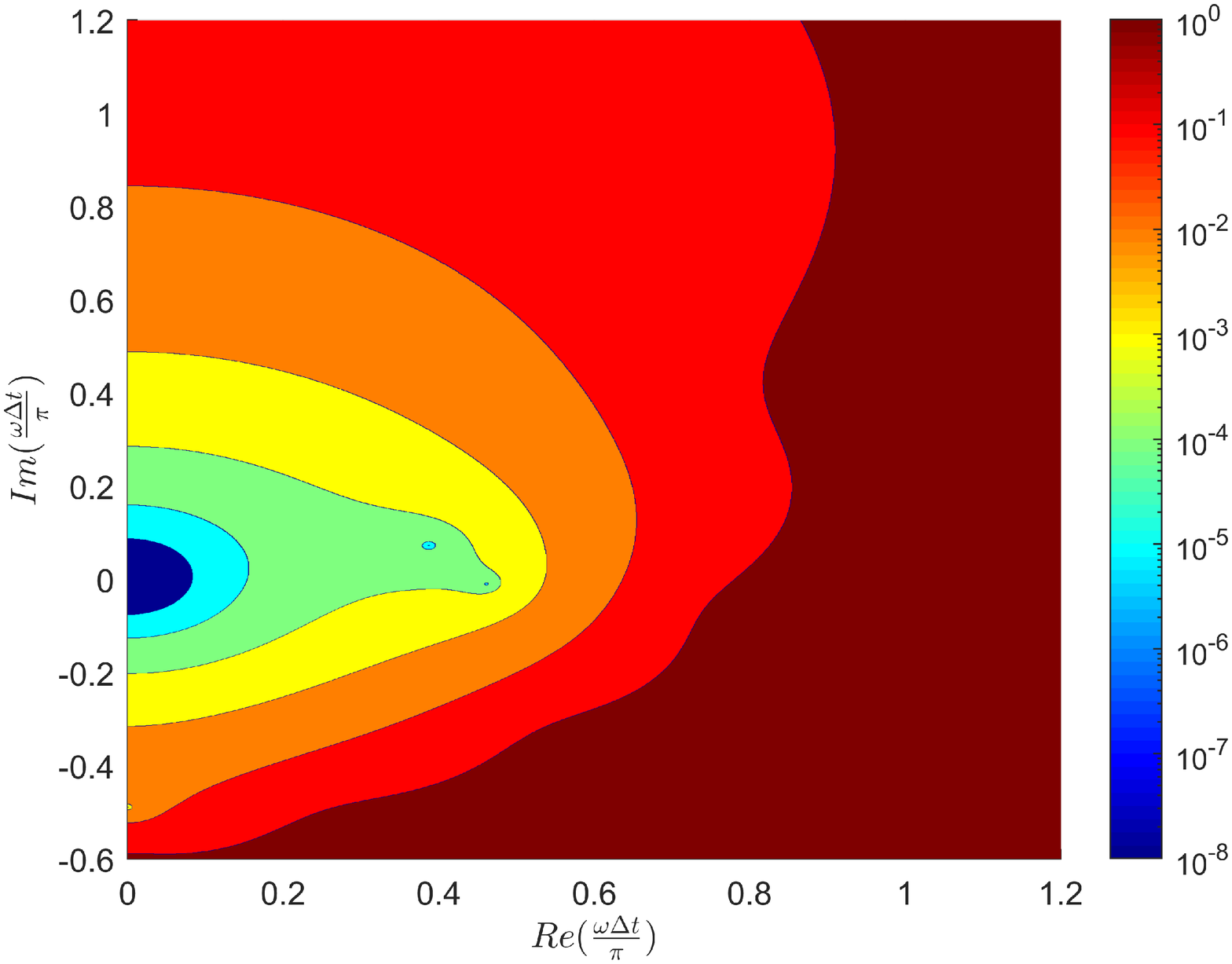}
    \caption{4th Order LDDRK46($\Delta t$)}
    \label{fig:LDDRK46v3}
  \end{subfigure}
  \hfill
  \begin{subfigure}[b]{0.49\textwidth}
    \includegraphics[width=\textwidth]{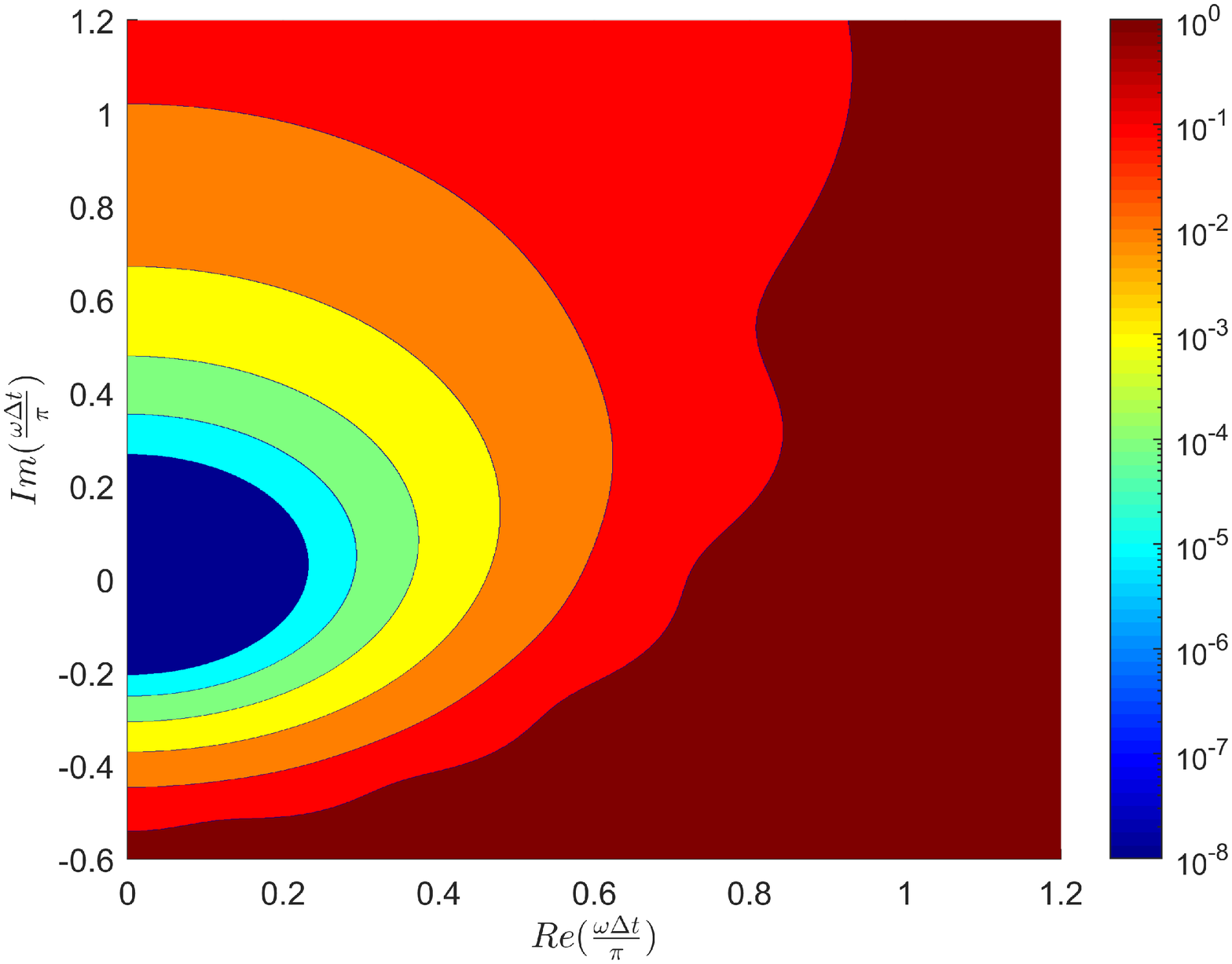}
    \caption{Maximal Order RK10($2\Delta t$)}
    \label{fig:twoTimeStepRK10}
  \end{subfigure}
  
  \begin{subfigure}[b]{0.49\textwidth}
    \includegraphics[width=\textwidth]{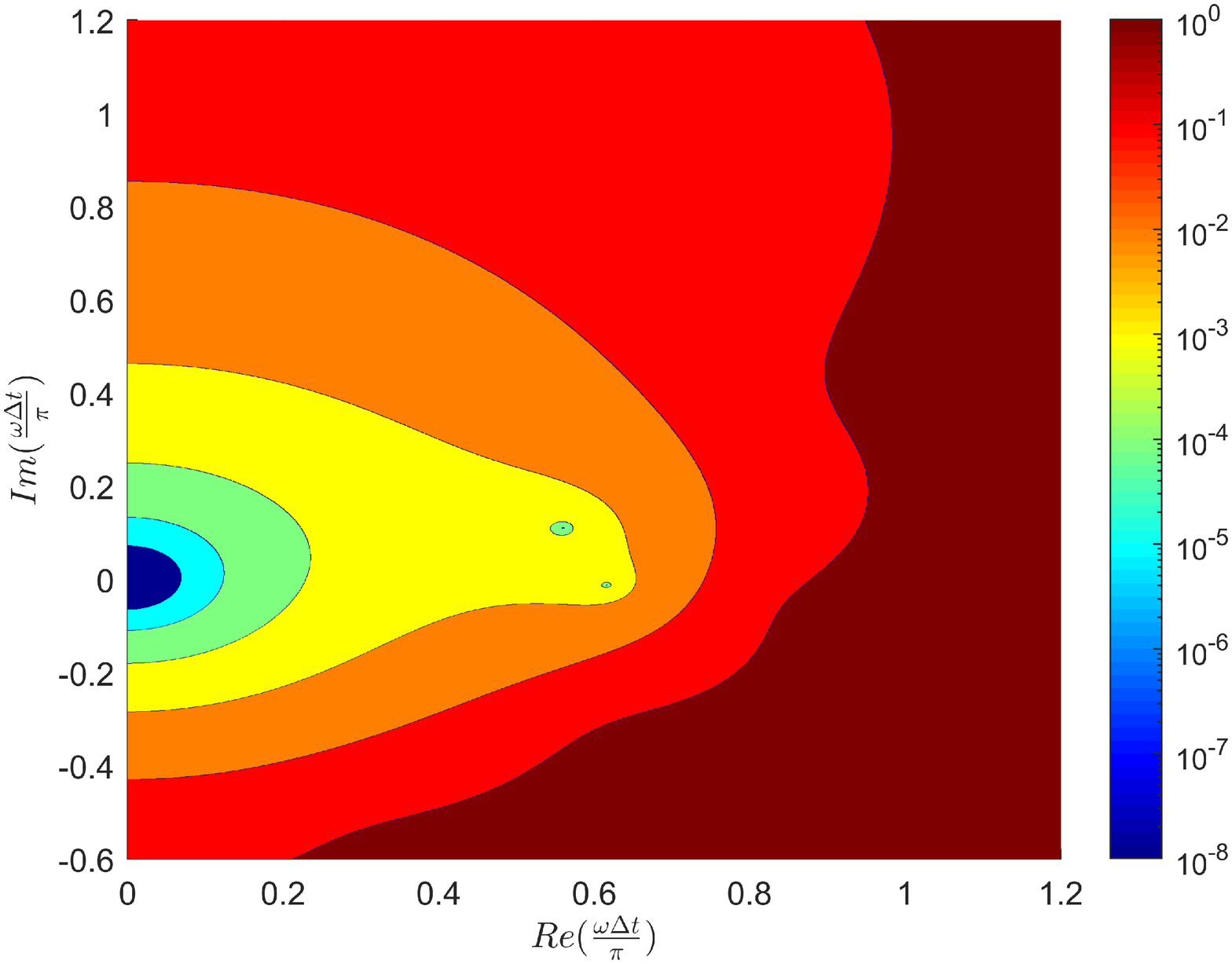}
    \caption{4th Order LDDRK56($\Delta t$)}
    \label{fig:LDDRK56v3}
  \end{subfigure}
  \hfill
  \begin{subfigure}[b]{0.49\textwidth}
    \includegraphics[width=\textwidth]{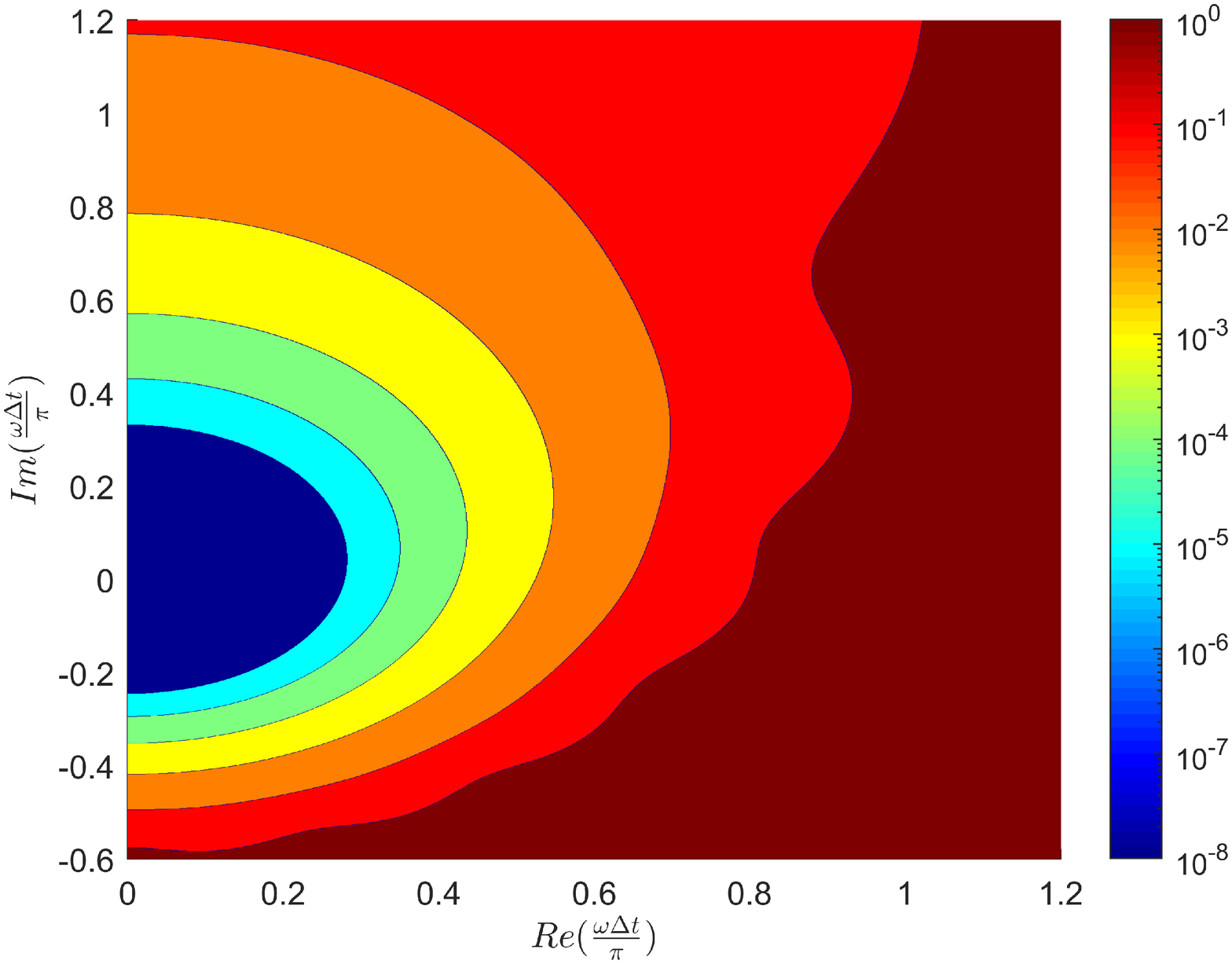}
    \caption{Maximal Order RK11($2\Delta t$)}
    \label{fig:twoTimeStepRK11}
  \end{subfigure}
  \caption{Plots of phase error $\eps_p$ in the complex $\odt$ plane for the LDDRK46 and LDDRK56 schemes of~\citet{hu+hussaini+manthey-1996}, and for maximal order 10-step and 11-step Runge--Kutta schemes with a time step of $2\Delta t$.}
  \label{fig:2StepLDDRKVStwoTimeStepRK}
\end{figure}
Figure~\ref{fig:2StepLDDRKVStwoTimeStepRK} shows the phase error for these schemes, $\eps_p(\odt)$, together with the comparable 10- and 11-stage maximal order schemes, $\eps_p(2\odt)$, while figure~\ref{fig:twoStepLDDRKvstwoTimeStepRkOnePercentPhaseError} compares which is the more accurate.
The maximal order schemes are clearly more accurate for the majority of values of $\odt$ despite using a twice as long time step as the optimized schemes.

\subsection{Stability considerations}

The phase error $\eps_p$ is not the only consideration when choosing a timestepping scheme.  For example, small values of $\eps_p$ are useless if the overall scheme is unstable.  Optimization of Runge--Kutta schemes are often restricted by a stability criterion.  For real $\omega$, $|r_e(\odt)| = \big|\exp\{-\I\odt\}\big| = 1$, meaning oscillation with no growth or decay.  A given scheme is said to be stable for $0 < \odt < \pi \eta_s$ if $|r(\odt)| < 1$ for $0< \odt<\pi \eta_s$; that is, the numerical scheme oscillates with possibly decay in amplitude, but with no growth.  In general, suppose a $p$-stage Runge--Kutta scheme is $q$-th order accurate with $q\geq 1$, so that
\begin{equation}
r(\odt) = 1 + \sum_{j=1}^q \frac{1}{j!}(-\I\odt)^j + \sum_{\mathclap{j=q+1}}^p c_j(-\I\odt)^j.
\end{equation}
For real $\omega$, in the limit $\odt\to0$, it is shown in~\ref{app:Stability} that
\begin{equation}
\Real\Big(\log\big(r(\odt)\big)\Big)  = \left\{\begin{array}{ll}\displaystyle
  (-1)^{n+1}\bigg[\!\left(c_{2n+2}-\frac{1}{(2n+2)!}\right)
  \\\displaystyle\qquad
  - \left(c_{2n+1}-\frac{1}{(2n+1)!}\right)\!\bigg](\odt)^{2n+2} + O\big((\odt)^{2n+4}\big)
  & \text{if $q=2n$,}\\\displaystyle
  (-1)^n\left(c_{2n}-\frac{1}{(2n)!}\right)(\odt)^{2n} + O\big((\odt)^{2n+2}\big)
  & \text{if $q=2n-1$.}
  \end{array}\right.
\label{equ:instability}
\end{equation}
Consequently, in order not to be unstable for arbitrarily small $\odt$, we require either $(-1)^{q/2}[(c_{q+1}-1/(q+1)!) - (c_{q+2}-1/(q+2)!)] < 0$ if $q$ is even or $(-1)^{(q+1)/2}(c_{q+1}-1/(q+1)!) < 0$ if $q$ is odd.  This is satisfied for all the optimized Runge--Kutta schemes of~\citet{bogey+bailly-2004}, and all the LDDRK schemes of~\citet{hu+hussaini+manthey-1996} apart from the LDDRK4 scheme (which is therefore slightly unstable for arbitrarily small real $\odt$).  For maximal order schemes, where $q=p$ and there is no flexibility in the choice of coefficients $c_j = 1/j!$, we find by setting $c_{q+1}=c_{q+2}=0$ in~\eqref{equ:instability} that they are stable for some $\eta_s>0$ if and only if $p=4m$ or $p=4m-1$ for some integer $m$.  In particular, this means that RK4, RK8, RK11, RK12 and RK16 are all stable for some $\eta_s$.

It is unclear what the equivalent restriction for complex frequencies should be.  One could consider the condition that $|r(\odt)/r_e(\odt)| < 1$, meaning that the numerical scheme gives a lower growth rate than the exact solution, although this is not necessarily a desirable property to have.  Here, we postulate that the notion of stability of a timestepping scheme is only relevant to real frequencies, while accuracy is important both for real and complex frequencies.

\subsection{Accuracy limits for real and complex frequencies}

By analogy with the real-frequency stability limit $\eta_s$ above, we may define the real-frequency accuracy limit $\eta_\delta$ such that $\eps_r < \delta$ for $0 < \odt < \pi\eta_\delta$.  Here, motivated by~\eqref{equ:accuracy} and a typical lower-bound order of magnitude $T/\Delta t \approx 100$, we will be particularly interested in $\delta = 10^{-4}$ and $\delta = 10^{-5}$, although we will also consider $\delta = 10^{-3}$ since this appears to be the error most optimized schemes have been optimized for.  Since we are also interested in the behaviour of timestepping schemes for non-constant-amplitude oscillations, corresponding to complex frequencies $\omega$, we define the analogous complex-frequency accuracy limit $\hat{\eta}_\delta$ such that $\eps_r < \delta$ for $0 < |\odt| < \pi\eta_\delta$; that is, an accuracy of $\eps_r < \delta$ is guaranteed irrespective of $\arg(\odt)$ provided the timestep is chosen sufficiently small that $|\odt| < \pi\eta_\delta$.

In order to make a sensible comparison between schemes with different numbers of stages, in what follows the amplification factor is re-scaled to be made comparable with RK4 such that the same number of stages will occur over the time domain $T$. For example, if a 12-stage scheme is given a timestep $3\Delta t$ then it will lead to a simulation with the same computational cost as a 4-stage scheme with a timestep $\Delta t$; in this case, the equivalent amplification factor for the 12-stage scheme will be $r(\odt) = \left(r_{12}(\omega 3 \Delta t)\right)^{1/3}$, so that after three $\Delta t$ ``time steps'' the effective amplification factor will be $r(\odt)^3 = r_{12}(\omega 3 \Delta t)$.  In general, the rescaled amplification factor for a $p-$stage Runge-Kutta scheme with unscaled amplification factor $r_p$ is
\begin{equation}
r(\odt) = \left ( r_p(\omega p \Delta t/4)\right )^{4/p};
\label{equ:rescaling}
\end{equation}
with the root minimising $\epsilon_r$ being chosen.  In order to minimize confusion, stability limits $\eta_s$ and accuracy limits $\eta_\delta$ and $\hat{\eta}_\delta$ calculated using the rescaled amplification factors will be denoted by $\lambda_s$, $\lambda_\delta$ and $\hat{\lambda}_\delta$ respectively; in order to make like-for-like comparisons between stability and accuracy limits for schemes with different numbers of stages, it is important to compare $\lambda_s$ and $\lambda_\delta$ rather than the unscaled $\eta_s$ and $\eta_\delta$ limits.

These accuracy limits $\lambda_\delta$ and $\hat{\lambda}_\delta$, along with the stability limits $\lambda_s$, are tabulated in table~\ref{table:rk}
\begin{table}
  \includegraphics[width=\linewidth, height=0.25\linewidth]{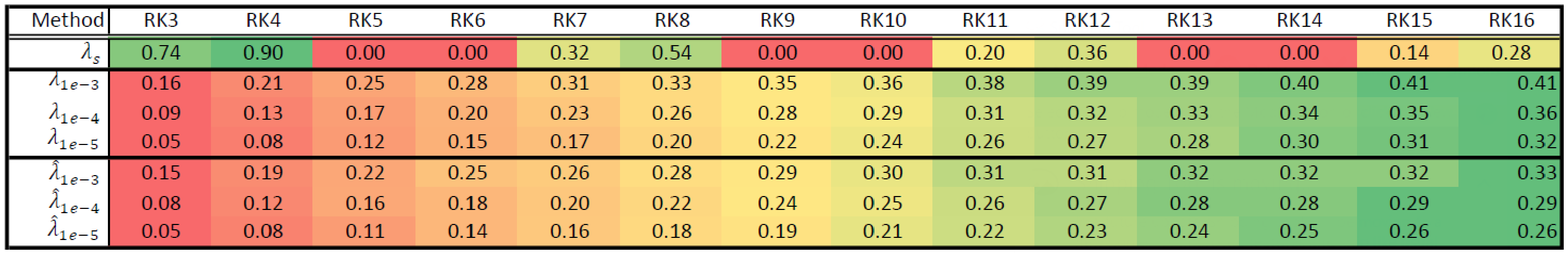}
\caption{Properties of rescaled amplification factors $r(\odt)$ from~\eqref{equ:rescaling}, meaning each scheme has the same computational cost, for $p$-stage maximal order Runge--Kutta timestepping schemes. Each cell is coloured by value, from highest (green) to lowest (red).}%
\label{table:rk}
\end{table}%
for maximal order Runge--Kutta schemes, and in table~\ref{table:lddrk}
\begin{table}
  \includegraphics[width=\linewidth]{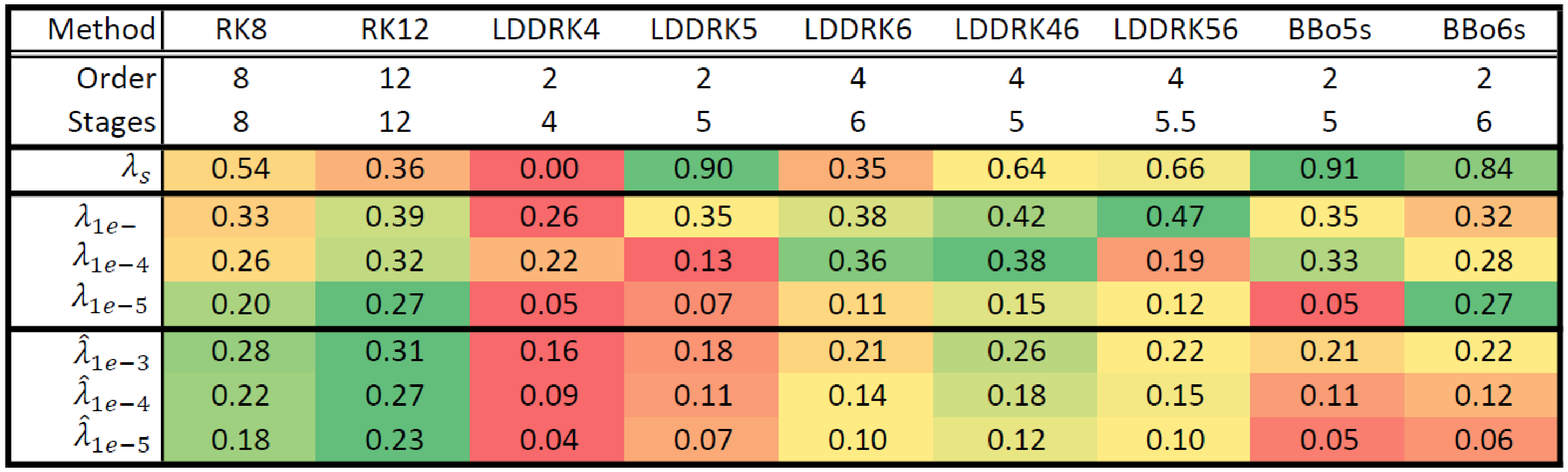}
  \caption{Properties of rescaled amplification factors $r(\odt)$ from~\eqref{equ:rescaling}, meaning each scheme has the same computational cost, for various optimized timestepping schemes.  Each cell is coloured by value, from highest (green) to lowest (red).}%
\label{table:lddrk}%
\end{table}%
for various optimized Runge--Kutta schemes.  For the maximal order schemes in table~\ref{table:rk}, lower order schemes are usually limited by accuracy, while higher order schemes are limited by their stability.  Behaviour for complex $\omega$ and for real $\omega$ are broadly comparable, although restricting to real $\omega$ does increase slightly the accuracy range. The situation is markedly different for the optimised schemes in table~\ref{table:lddrk}, where the optimised schemes generally outperform the maximal order schemes for the real-$\omega$ $\delta=10^{-3}$ case, and underperform their maximal order equivalents when either $\omega$ is complex, or when the desired error is $\delta = 10^{-5}$.  This verifies that the optimised schemes have been optimised solely for real $\omega$, with optimization parameters tuned to prioritize errors of around $10^{-3}$.

%%%%%%%%%%%%%%%%
%Reoptimization%
%%%%%%%%%%%%%%%%

\section{Optimisation for non-constant-amplitude oscillations}
\label{section:optimization}

All the optimisations performed to generate optimized schemes to date assume that $\omega$ is real. Therefore, there is an implicit assumption that the amplitude of the oscillations do not change with time. As discussed in section \ref{section:TheoreticalComparison}, the fact that the considered optimised schemes fail to outperform their maximal order counterparts except for very small regions in the complex plane centered about the real axis is the manifestation of this implicit assumption.  In this section, we consider optimizations over complex $\omega$.

\subsection{Optimisation metrics for a range of complex \texorpdfstring{$\omega$}{ω}}

\begin{figure}
  \begin{subfigure}[b]{0.49\textwidth}
     \centering
     \psfrag{a1}[r][r]{$\alpha_1\eta$}
     \psfrag{a2}[r][r]{$\alpha_2\eta$}
     \psfrag{eta}[lB][lB]{$\eta$}
     \includegraphics[height=0.85\textwidth]{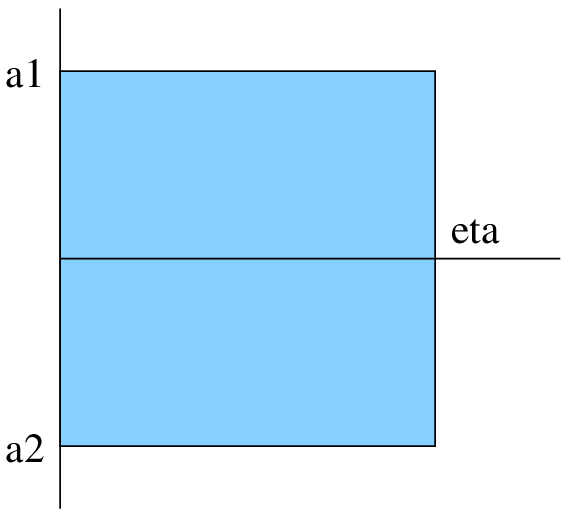}
\caption{}
    \label{fig:SquareGeometry}
  \end{subfigure}
  \hfill
  \begin{subfigure}[b]{0.49\textwidth}
     \centering
     \psfrag{b1}[lB][lB]{$\beta_1$}
     \psfrag{b2}[lt][lt]{$\beta_2$}
     \psfrag{eta}[lB][lB]{$\eta$}
     \includegraphics[height=0.85\textwidth]{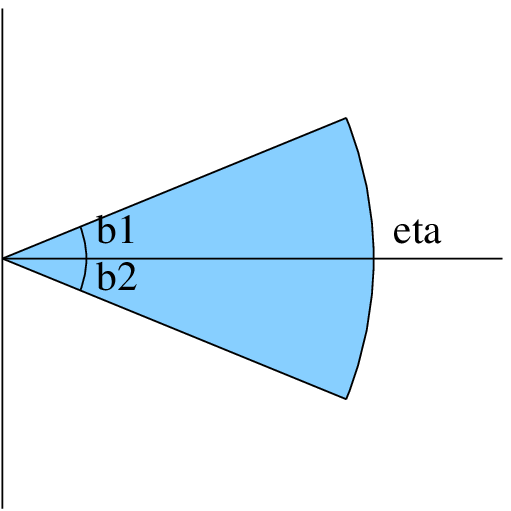}
    \caption{}
    \label{fig:SectorGeometry}
  \end{subfigure}
  \caption{Schematic of the two regions of the complex-$\odt/\pi$ plane over which the optimization integrals~\eqref{equ:metricSquare} and~\eqref{equ:metricSector} are performed.}
  \label{figure:geometry}
\end{figure}
Inspired by previous such optimizations for spatial derivatives~\citep{brambley+markeviciute-2017-aiaa}, we consider two generalizations of the optimization metric~(\ref{equ:metric}) to regions of the complex $\odt$ plane, shown schematically in figure~\ref{figure:geometry}.  For the rectangular region in figure~\ref{fig:SquareGeometry}, defined by $\alpha_1 \geq 0, \alpha_2 \leq 0$ and $\eta > 0$, the optimization is taken to be 
\begin{subequations}\begin{align}
  e &= \frac{1}{(|\alpha_1|+ |\alpha_2|)\pi\eta} \int_{\alpha_2\pi\eta} ^{\alpha_1\pi\eta} \int_0^{\pi\eta}\left|r(p + \I q) - r_e(p + \I q)) \right|^2\, \intd p \intd q,\\
  E &= \frac{1}{(|\alpha_1|+ |\alpha_2|)\pi\eta} \int_{\alpha_2\pi\eta} ^{\alpha_1\pi\eta} \int_0^{\pi\eta}\Big|\bar{\omega}\big((p + \I q)/\Delta t\big)\Delta t - (p+\I q) \Big|^2\, \intd p \intd q.
\end{align}\label{equ:metricSquare}\end{subequations}
The choice of $\alpha_1$ and $\alpha_2$ allows the creation of optimised schemes specialised for growth and decay rates \emph{per timestep} between $\exp{\{\alpha_1 \pi\eta\}}$ and $\exp{\{\alpha_2 \pi\eta\}}$. As $\alpha_1,\alpha_2\to0,$ the original metric~\eqref{equ:metric} is recovered. The second geometry considered is a sector with radius $\pi\eta$ and angle $\left |\beta_1\right | + \left |\beta_2 \right |$. Letting $\odt$ be of the form $\rho\exp{\{\I\theta\}}$ gives the optimization metric as
\begin{subequations}\begin{align}
    e &= \frac{1}{(|\beta_1| + |\beta_2|)\pi\eta} \int_{\beta_2} ^{\beta_1} \int_0^{\pi\eta}|r(\rho\e^{\I\theta}) - r_e(\rho\e^{\I\theta})|^2\, \rho\intd\rho\intd\theta \label{equ:sectormetric}\\
    E &= \frac{1}{(|\beta_1| + |\beta_2|)\pi\eta} \int_{\beta_2} ^{\beta_1} \int_0^{\pi\eta}|\bar{\omega}(\rho\e^{\I\theta}\!/\Delta t)\Delta t - \rho\e^{\I\theta}|^2\, \rho\intd\rho\intd\theta
\end{align}\label{equ:metricSector}\end{subequations}
The choice of $\beta_1$ and $\beta_2$ allows the creation of optimised schemes specialised for growth and decay rates \emph{per period} of between $\exp{\{2\pi\tan\beta_1 \}}$ and $\exp{\{2\pi\tan\beta_2 \}}$.  As $\beta_1,\beta_2\to0$, a variation of metric (\ref{equ:metric}) is obtained which penalises larger errors away from the origin.

\subsection{Results of complex-\texorpdfstring{$\omega$}{ω} optimizations}

The metrics described above were used to create various p-stage optimised Runge-Kutta schemes of order 4. The optimiser was constrained to ensure stability by requiring that  
\begin{align}
    & (c_5 - \frac{1}{5!}) - (c_6 - \frac{1}{6!}) < 0,&
  &\text{and}&
    & \eta_s > \hat{\eta}_s&
\end{align}
for some given $\hat{\eta}_s$, with the first constraint taken from~\eqref{equ:instability}.

Many optimization parameters were considered, and only three of the better performing schemes are described here.  These schemes, a 6-, 8- and 12-stage scheme labelled Opt6, Opt8 and Opt12 here, are fairly typical of the types of optimized schemes obtained.  All were obtained using the sector optimization~\eqref{equ:sectormetric}, although similar results were also obtained from the other metrics.  The parameters used for the optimization, and the coefficients obtained, are given in table~\ref{table:OptimizedCoefficients}.
\begin{table}%
\centering%
\setlength{\tabcolsep}{22pt}%
\begin{tabular}{ c c c c}
 & Opt6 & Opt8 & Opt12 \\ 
 \hline \hline
 $\eta$ & $1/2$ & $3/4$ & $1$\\
 $\beta_1$ &$\pi/6$&$\pi/6$&$\pi/6$ \\
 $\beta_2$ &$-\pi/6$&$-\pi/6$&$0$ \\
 $\hat{\eta}_s$ &$1/2$&$1$&$1/2$ \\
\hline
 c$_5$ & $7.86006019\times 10^{-3}$	 & $8.27554045\times 10^{-3}$ & $8.33315438\times 10^{-3}$	 \\  
 c$_6$ & $1.21477435\times 10^{-3}$	 & $1.37185292\times 10^{-3}$	 & $1.38885733\times 10^{-3}$	 \\  
 c$_7$ &  & $1.76272985\times 10^{-4}$	 & $1.98395863\times 10^{-4}$	 \\  
 c$_8$ &  & $2.05839623\times 10^{-5}$	  & $2.47338621\times 10^{-5}$	 \\  
 c$_9$ &  &  & $2.75123146\times 10^{-6}$	 \\  
 c$_{10}$ &  &  & $2.65593613\times 10^{-7}$	 \\  
 c$_{11}$ &  &  & $2.28460890\times 10^{-8}$	 \\  
 c$_{12}$ &  &  & $1.65356900\times 10^{-9}$  \\
 \hline
\end{tabular}%
\caption{Details of the optimization parameters used and coefficients obtained for the Opt6, Opt8 and Opt12 4th order Runge--Kutta schemes.  Note that, for $1\leq j\leq 4$, 4th order accuracy requires that $c_j = 1/j!$.}
\label{table:OptimizedCoefficients}
\end{table}

\begin{figure}
\makebox[\linewidth][c]{
  \begin{subfigure}[b]{0.3\textwidth}
  \centering
    \includegraphics[width=\textwidth]{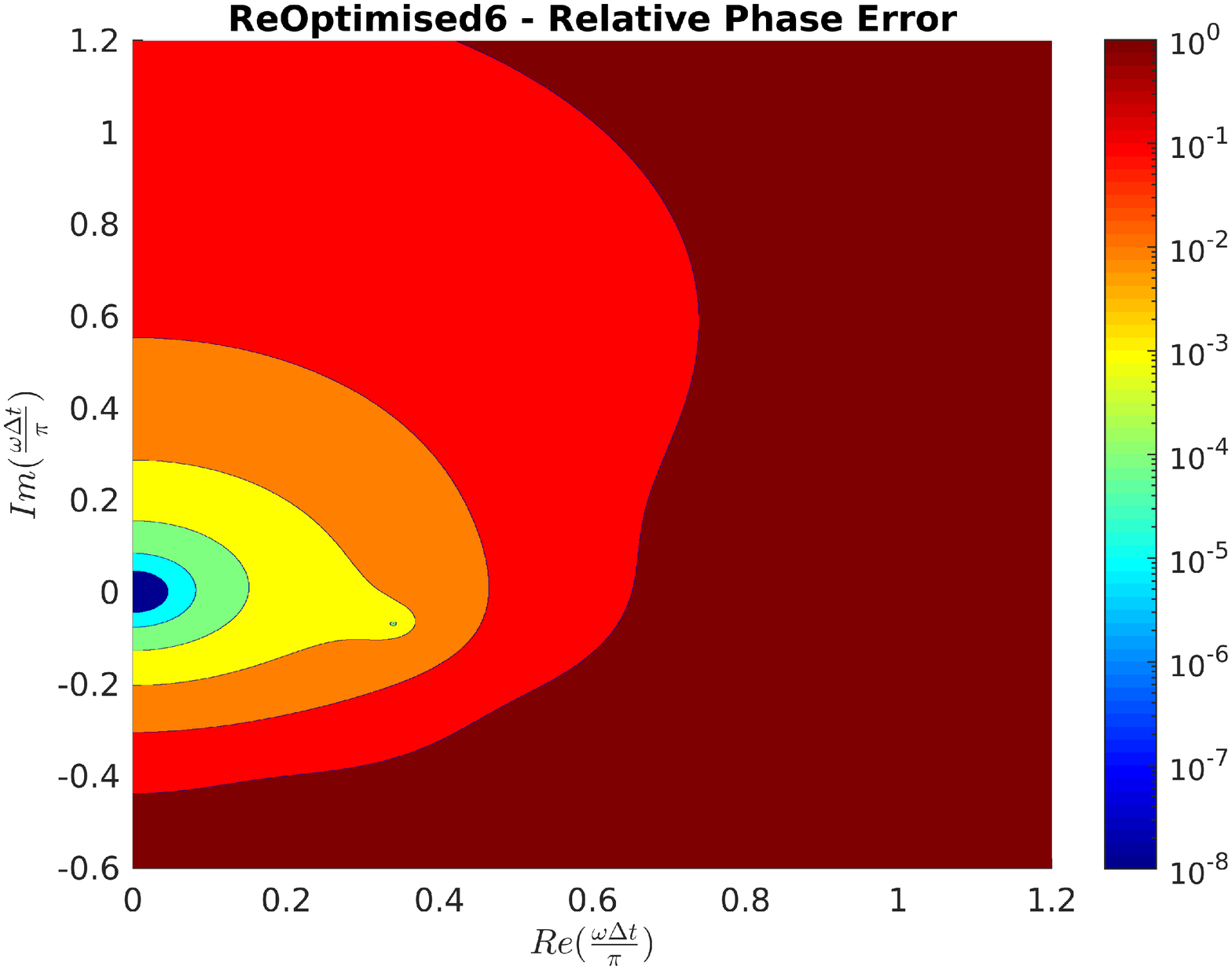}
    \caption{\small% $\beta_1 = \frac{\pi}{6}$ $\beta_2 = \frac{-\pi}{6}$, $\eta_s = \frac{\pi}{2}
    Opt6}
    \label{fig:Reoptimised6StagePhaseError}
  \end{subfigure}
  \hfill
  \begin{subfigure}[b]{0.3\textwidth}
  \centering
    \includegraphics[width=\textwidth]{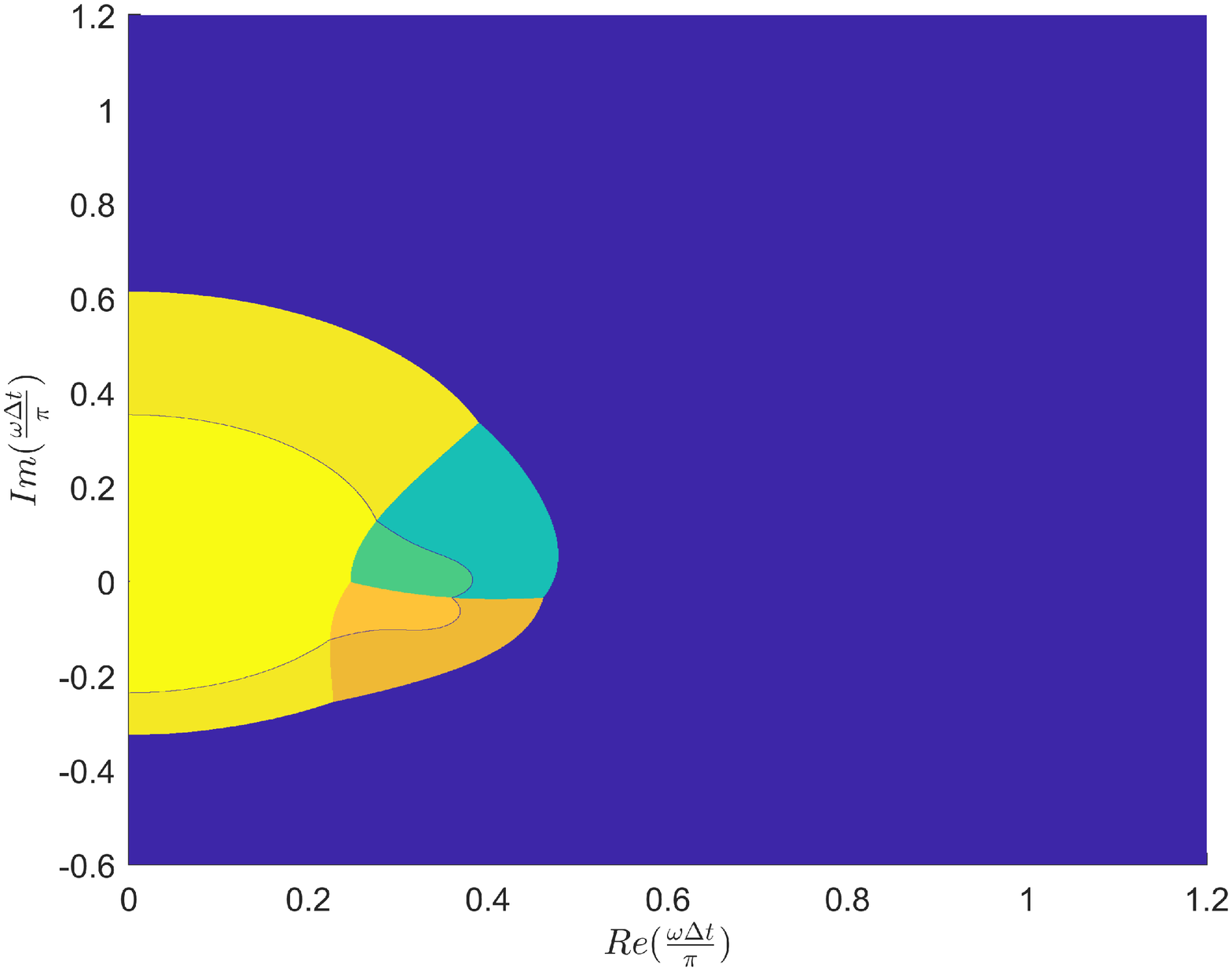}
    \caption{\small \textcolor{yellow}{RK6} vs \textcolor{green}{LDDRK6} vs \textcolor{orange}{Opt6}}
    \label{fig:Reoptimised6StagePhaseErrorComparison}
  \end{subfigure}
  \hfill
   \begin{subfigure}[b]{0.3\textwidth}
  \centering
    \includegraphics[width=\textwidth]{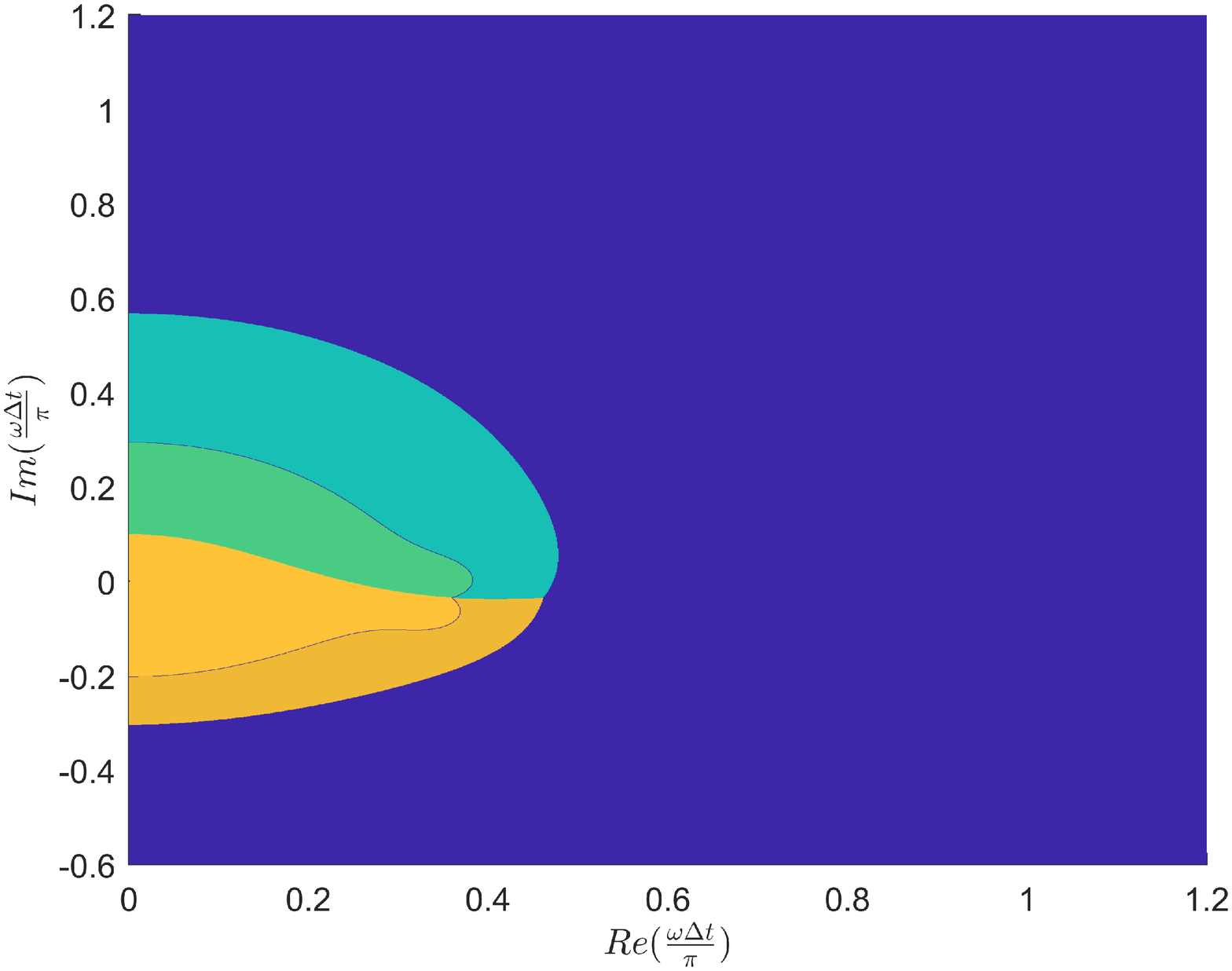}
    \caption{\small \textcolor{green}{LDDRK6} vs \textcolor{orange}{Opt6}}
    \label{fig:Reoptimised6StagePhaseErrorComparisonVSLDDRK6Only}
  \end{subfigure}
  }
 \makebox[\linewidth][c]{
  \begin{subfigure}[b]{0.3\textwidth}
  \centering
    \includegraphics[width=\textwidth]{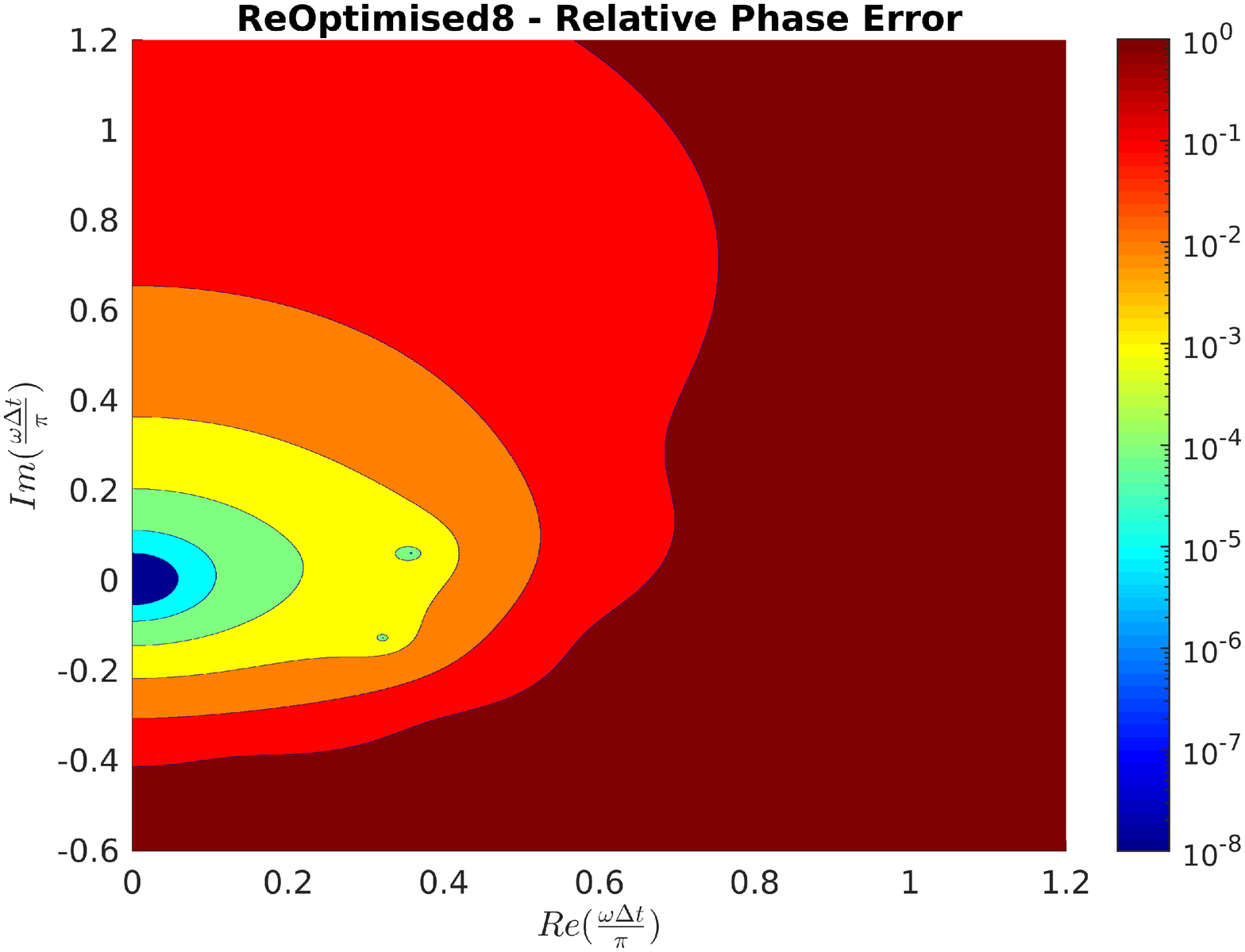}
    \caption{\small% $\beta_1 = \frac{\pi}{6}$ $\beta_2 = \frac{-\pi}{6}$, $\eta_s = \frac{\pi}{1}$
    Opt8}
    \label{fig:Reoptimised8StagePhaseError}
  \end{subfigure}
  \hfill
  \begin{subfigure}[b]{0.3\textwidth}
  \centering
    \includegraphics[width=\textwidth]{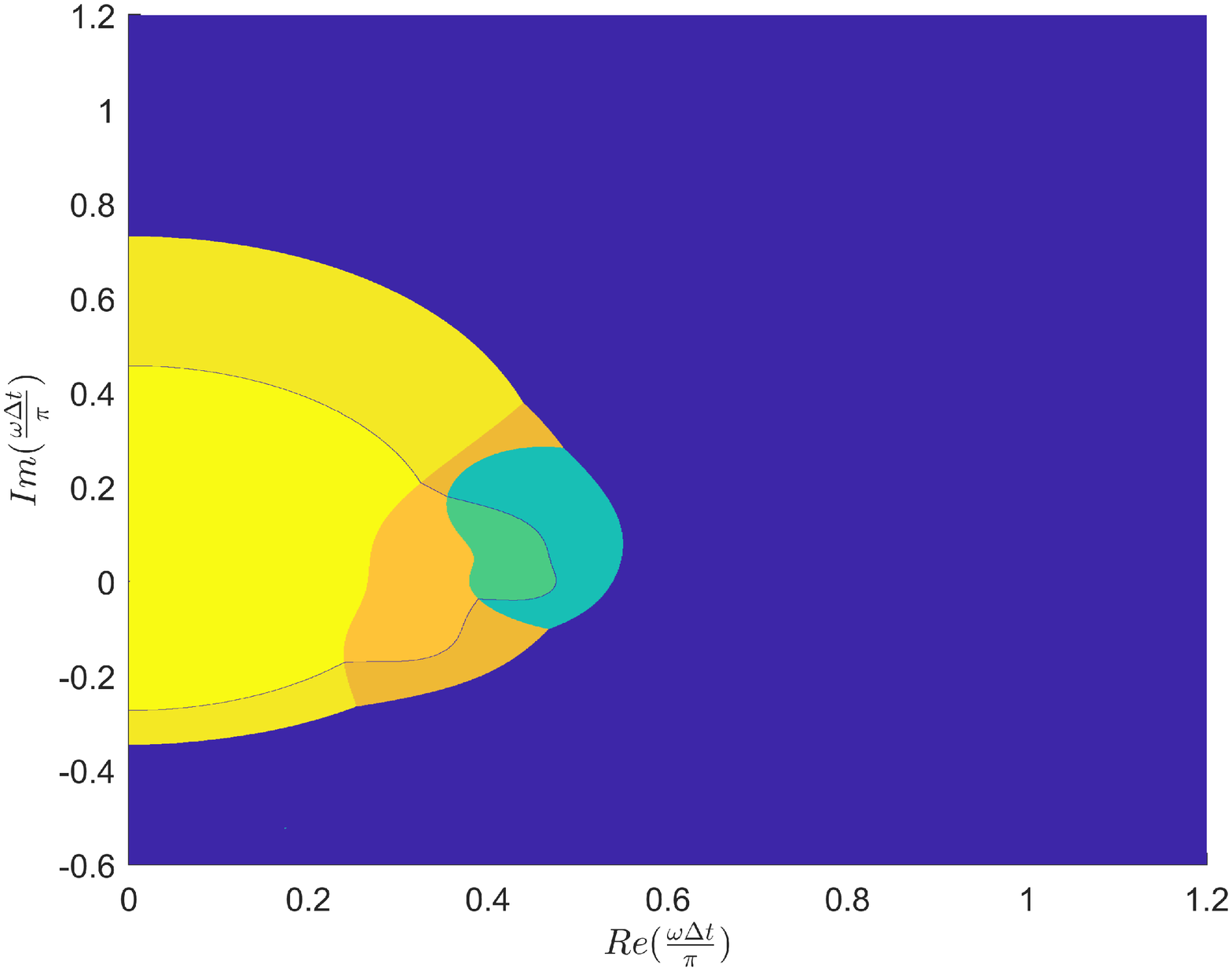}
    \caption{\small \textcolor{yellow}{RK8} vs \textcolor{green}{LDDRK56} vs \textcolor{orange}{Opt8}}
    \label{fig:Reoptimised8StagePhaseErrorComparison}
  \end{subfigure}
  \hfill
  \begin{subfigure}[b]{0.3\textwidth}
  \centering
    \includegraphics[width=\textwidth]{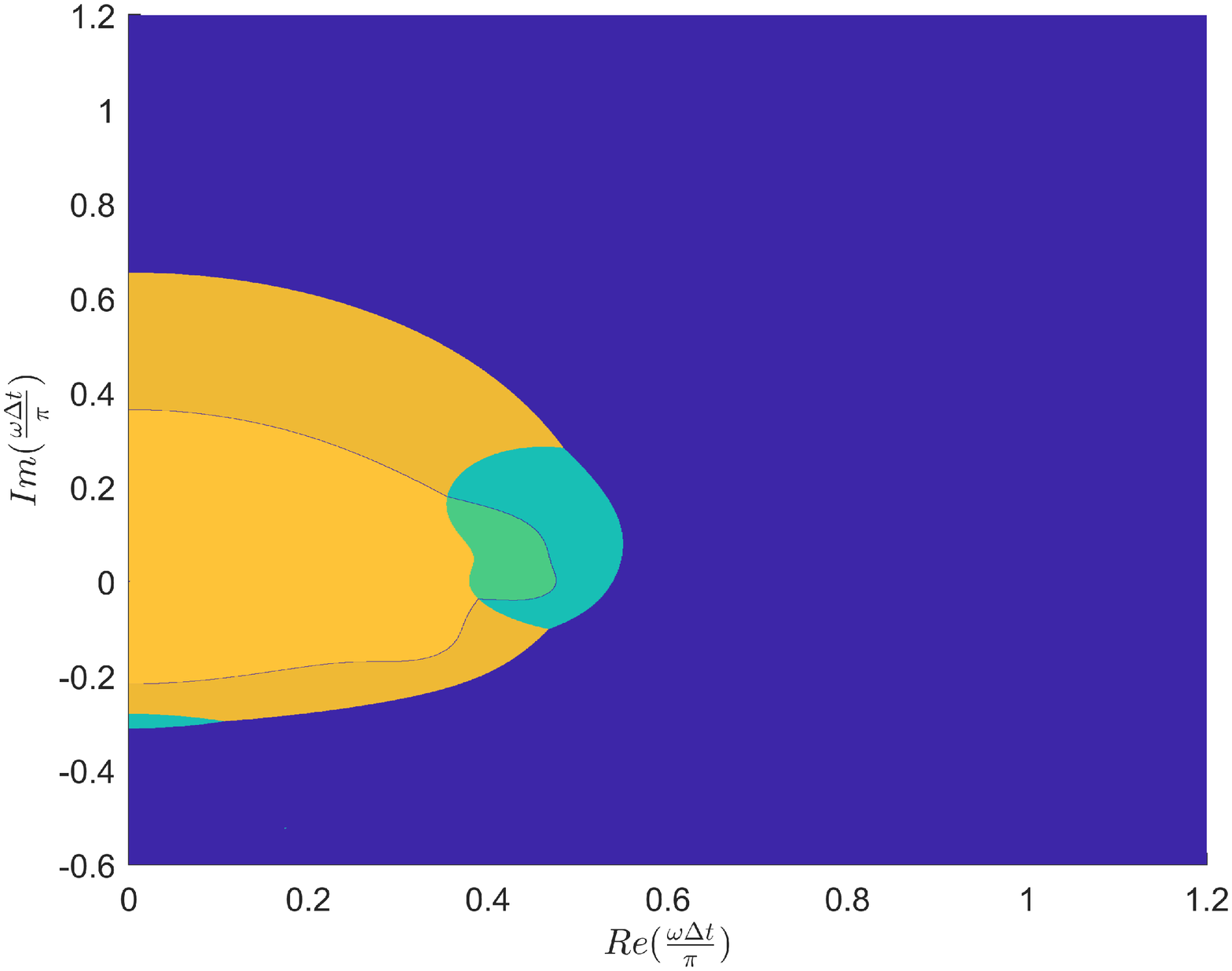}
    \caption{\small \textcolor{green}{LDDRK56} vs \textcolor{orange}{Opt8}}
    \label{fig:Reoptimised8StagePhaseErrorComparisonVSLDDRK56Only}
  \end{subfigure}
  }
  \makebox[\linewidth][c]{
  \begin{subfigure}[b]{0.3\textwidth}
  \centering
    \includegraphics[width=\textwidth]{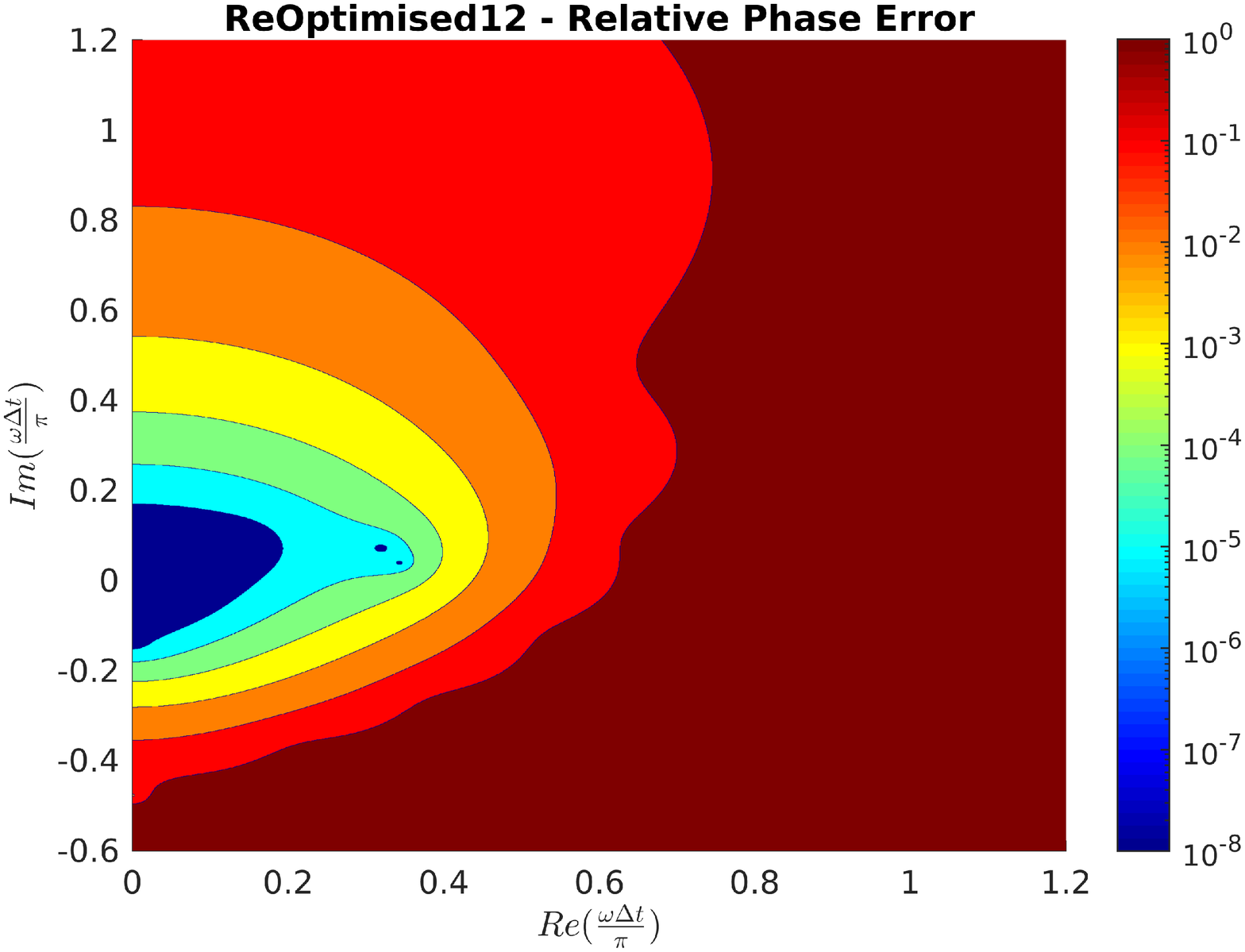}
    \caption{\small% $\beta_1 = \frac{\pi}{6}$ $\beta_2 = \frac{-\pi}{6}$, $\eta_s = \frac{\pi}{4}$
    Opt12}
    \label{fig:Reoptimised12StagePhaseError}
  \end{subfigure}
  \hfill
  \begin{subfigure}[b]{0.3\textwidth}
  \centering
    \includegraphics[width=\textwidth]{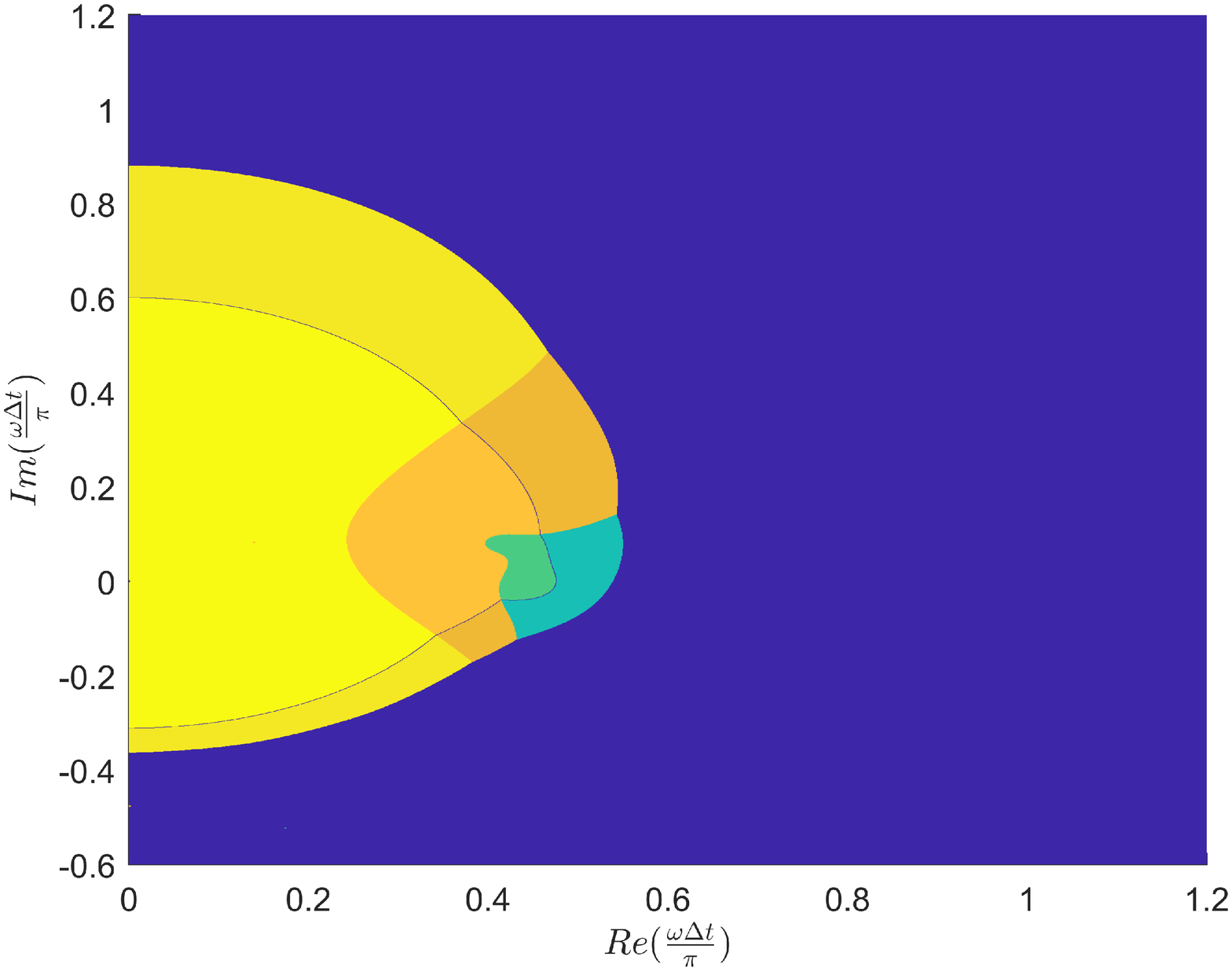}
    \caption{\small \textcolor{yellow}{RK12} vs \textcolor{green}{LDDRK56} vs \textcolor{orange}{Opt12}}
    \label{fig:Reoptimised12StagePhaseErrorComparison}
  \end{subfigure}
  \hfill
  \begin{subfigure}[b]{0.3\textwidth}
  \centering
    \includegraphics[width=\textwidth]{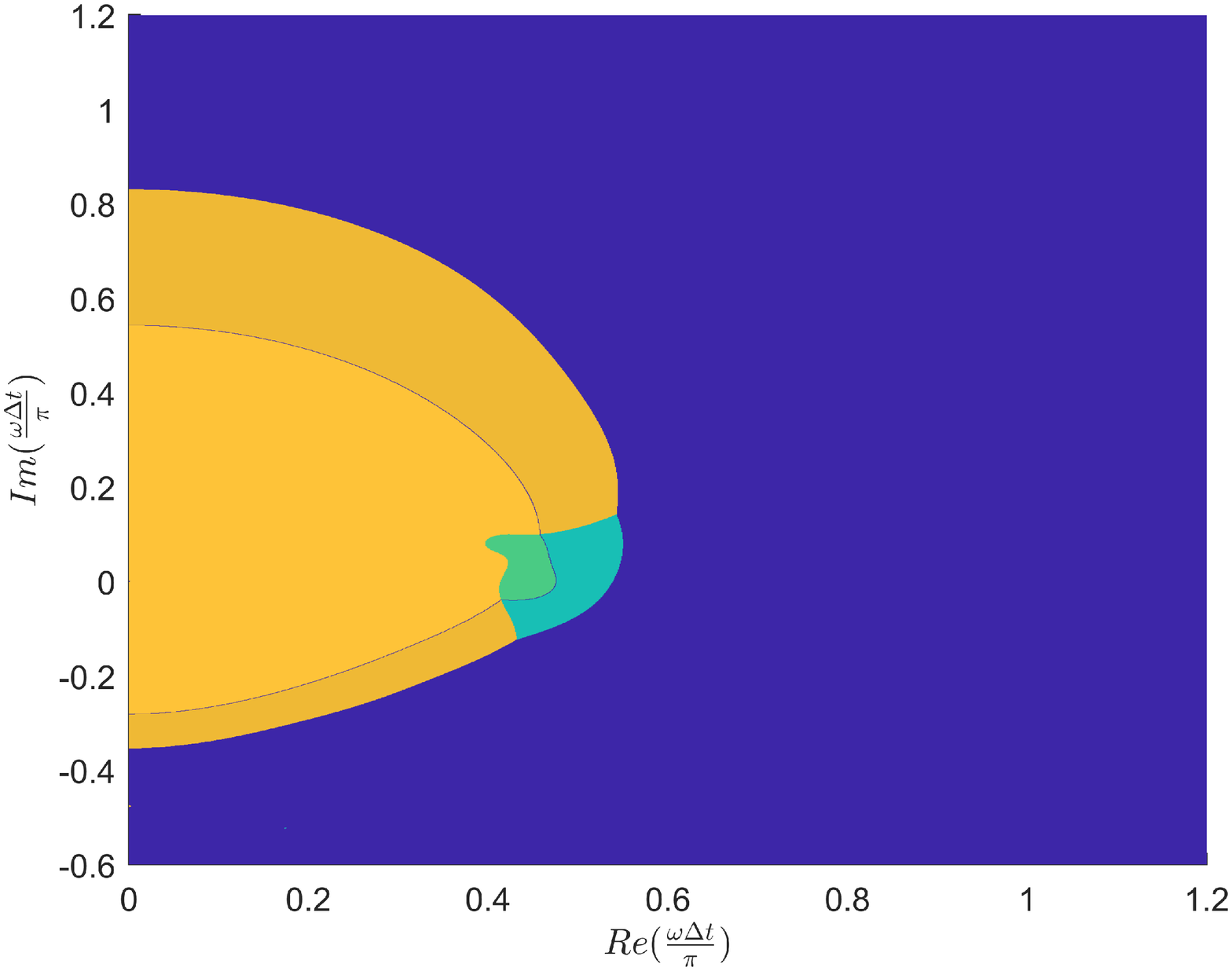}
    \caption{\small \textcolor{green}{LDDRK56} vs \textcolor{orange}{Opt12}}
    \label{fig:Reoptimised12StagePhaseErrorComparisonVSLDDRK56Only}
  \end{subfigure}
  }
  \caption{\small Plots of $\eps_p(\odt p/4)$ for $p$-stage optimized schemes.  The rescaling of $\odt$ by $p/4$ means plots for different numbers of stages are directly comparable. (\subref{fig:Reoptimised6StagePhaseError}),(\subref{fig:Reoptimised8StagePhaseError}) and (\subref{fig:Reoptimised12StagePhaseError})  are plots of $\eps_p$ on a logarithmic scale. (\subref{fig:Reoptimised6StagePhaseErrorComparison}),(\subref{fig:Reoptimised6StagePhaseErrorComparisonVSLDDRK6Only}),(\subref{fig:Reoptimised8StagePhaseErrorComparison}),(\subref{fig:Reoptimised8StagePhaseErrorComparisonVSLDDRK56Only}),(\subref{fig:Reoptimised12StagePhaseErrorComparison}),(\subref{fig:Reoptimised12StagePhaseErrorComparisonVSLDDRK56Only}) are phase error comparison plots where the colour indicates the scheme with the lowest phase error.  Lighter colours indicate an error of less than $0.1\%$, while violet indicates no scheme is $1\%$ accurate.}
    \label{fig:ReoptimisedPhaseErrors}
\end{figure}%
Figure \ref{fig:ReoptimisedPhaseErrors} shows the phase errors, $\epsilon_p$, of various optimised and maximal order Runge-Kutta schemes.  The phase errors are plotted such that the amplification factors have been rescaled to be comparable with RK4, as given in~\eqref{equ:rescaling}.  As can be seen, the newly optimized schemes outperform their LDDRK counterparts for both non-constant and constant amplitude oscillations.  However, they still fail to perform better than their maximal order counterparts for a range of complex values of $\omega$.  It may also be noted that the stability constraint tended to be a limiting factor for the constrained optimiser;  if the stability requirements were too large, the optimiser failed to minimise the integral or produce a scheme that offered any performance improvements.

\begin{table}
  \includegraphics[width=\linewidth]{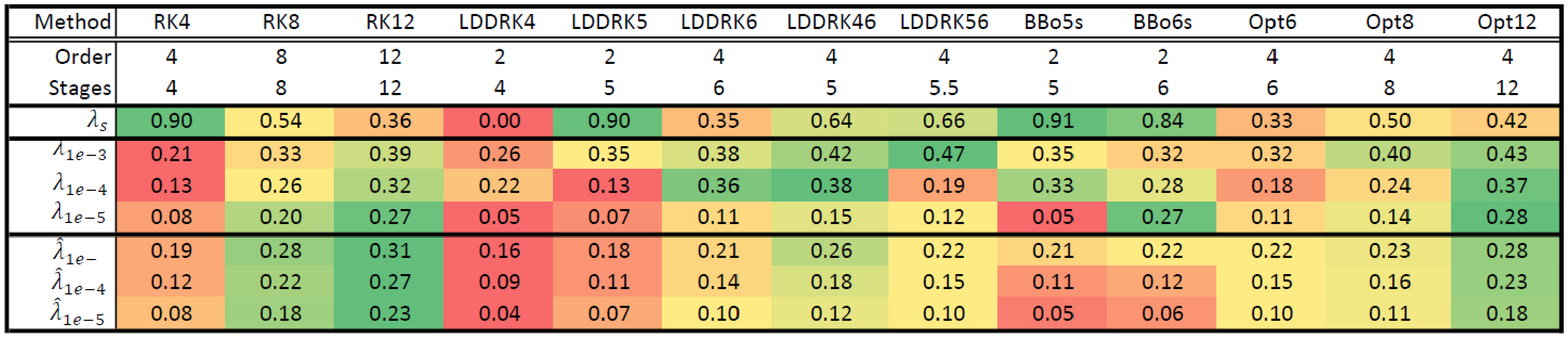}
  \caption{Properties of rescaled amplification factors $r(\odt)$ from~\eqref{equ:rescaling} . Each cell is coloured by value, from highest (green) to lowest (red).}%
  \label{tbl:StabilityReoptimisedSchemes}
\end{table}
Table~\ref{tbl:StabilityReoptimisedSchemes} shows the re-scaled stability limits for various maximal order and optimised schemes; because of the rescaling, all schemes have the same computational cost. The table clearly shows that lower order methods, except for the re-optimised 6 stage method, have the highest stability whilst simultaneously having the lowest accuracy limits comparatively. LLDRK5, LDDRK6, LDDRK46 and reoptimized 8 and 12 stage methods have similar, in some cases larger stability and $0.1\%$ accuracy limits than the maximal order 12-stage Runge-Kutta scheme. It may also be noted that only the 12-stage re-optimised scheme offers similar accuracy limits for complex $\omega$ to the 8th and 12th maximal order schemes.

%%%%%%%%%%%%%%%%
%Test Case%
%%%%%%%%%%%%%%%%

\section{Comparison using a realistic 1D test case}
\label{section:1d}

We now investigate how the theoretical behaviour described above translates into performance in practice, by comparing the performance of the various timestepping schemes for the simple 1D wave propagation problem from Ref.~\citenum{brambley-2016-jcp}.  The problem to be solved is 
\begin{align}
\frac{\partial p}{\partial t} + \frac{\partial v}{\partial x} &= -k_p(x) p, &
\frac{\partial v}{\partial t} + \frac{\partial p}{\partial x} &= -k_v(x) v.
\label{equ:problem}
\end{align}
These equations support wave propagation in both the
positive and negative $x$-directions at a wave speed of $1$.   Equation~\eqref{equ:problem} is solved on a periodic $x$-domain $[0,\,24)$, with initial conditions $v(x,0) = p(x,0)$ and damping $k_p(x) = k_v(x)$ as specified in Ref.~\citenum{brambley-2016-jcp} consisting of a wave packet with wavelength $1$ propagating across a damping region of length $2$ and decaying by a factor of $\e^{-6}$.  By comparing with the analytic solution $p_a(x,t)$, $v_a(x,t)$, given in~\citep{brambley-2016-jcp}, the numerical error is then given by
\begin{align}
  \mathrm{Error} &= \dfrac{\sup_{x\in[0,24)}\Big\{\big|p(x,T) - p_a(x,T)\big|,\,\big|v(x,T) - v_a(x,T)\big|\Big\}}{\sup_{x\in[0,24)}\Big\{\big|p_a(x,T)\big|,\,\big|v_a(x,T)\big|\Big\}} &
&\text{with }T=24.      
\label{equ:error}      
\end{align}
  
\begin{figure}%
\centering%
\psfrag{Error}[B][B]{Error}%  
\psfrag{dt}[t][t]{$\Delta t$}%  
\psfrag{CFL}[B][B]{CFL}%  
\psfrag{Effort}[t][t]{Effort}%
\psfrag{y1.0e-12}[r][r]{$10^{-12}$}%  
\psfrag{y1.0e-09}[r][r]{$10^{-9}$}%  
\psfrag{y1.0e-06}[r][r]{$10^{-6}$}%  
\psfrag{y1.0e-03}[r][r]{$10^{-3}$}%  
\psfrag{y1.0e-02}[r][r]{$10^{-2}$}%  
\psfrag{y1.0e-01}[r][r]{$10^{-1}$}%  
\psfrag{y2.0e-01}[r][r]{$0.2$}%  
\psfrag{y5.0e-01}[r][r]{$0.5$}%  
\psfrag{y1.0e+00}[r][r]{$1$}%  
\psfrag{xt1.0e-02}[B][B]{$10^{-2}$}%  
\psfrag{xt5.0e-02}[B][B]{$0.05$}%  
\psfrag{xt1.0e-01}[B][B]{$10^{-1}$}%  
\psfrag{xt2.0e-01}[B][B]{$0.2$}%  
\psfrag{xt5.0e-01}[B][B]{$0.5$}%  
\psfrag{xt1.0e+00}[B][B]{$1$}%  
\psfrag{xt2.0e+00}[B][B]{$2$}%  
\psfrag{xt4.0e+00}[B][B]{$4$}%  
\psfrag{x1.0e+09}[B][B]{$10^9$}%  
\psfrag{x1.0e+08}[B][B]{$10^8$}%  
\psfrag{x1.0e+07}[B][B]{$10^7$}%  
\psfrag{x1.0e+06}[B][B]{$10^6$}%  
\psfrag{x5.0e-01}[B][B]{$0.5$}%  
\psfrag{x2.0e-01}[B][B]{$0.2$}%  
\psfrag{x1.0e-01}[B][B]{$10^{-1}$}%  
\psfrag{x5.0e-02}[B][B]{$5\times 10^{-2}$}%  
\psfrag{x2.0e-02}[B][B]{$2\times 10^{-2}$}%  
\psfrag{x1.0e-02}[B][B]{$10^{-2}$}%  
\psfrag{x5.0e-03}[B][B]{$5\times 10^{-3}$}%  
\psfrag{x1.0e-03}[B][B]{$10^{-3}$}%  
\psfrag{RK4}[r][r]{RK4}%  
\psfrag{RK8}[r][r]{RK8}%  
\psfrag{RK11}[r][r]{RK11}%  
\psfrag{RK12}[r][r]{RK12}%  
\psfrag{RK16}[r][r]{RK16}%  
\psfrag{BBo5s}[r][r]{BBo5s}%  
\psfrag{BBo6s}[r][r]{BBo6s}%  
\psfrag{PET6}[r][r]{Opt6}%  
\psfrag{PET8}[r][r]{Opt8}%  
\psfrag{PET12}[r][r]{Opt12}%  
\psfrag{LDDRK4}[r][r]{\small{LDDRK4}}%  
\psfrag{LDDRK5}[r][r]{\small{LDDRK5}}%  
\psfrag{LDDRK6}[r][r]{\small{LDDRK6}}%  
\psfrag{LDDRK46}[r][r]{\small{LDDRK46}}%  
\psfrag{LDDRK56}[r][r]{\small{LDDRK56}}%  
\includegraphics{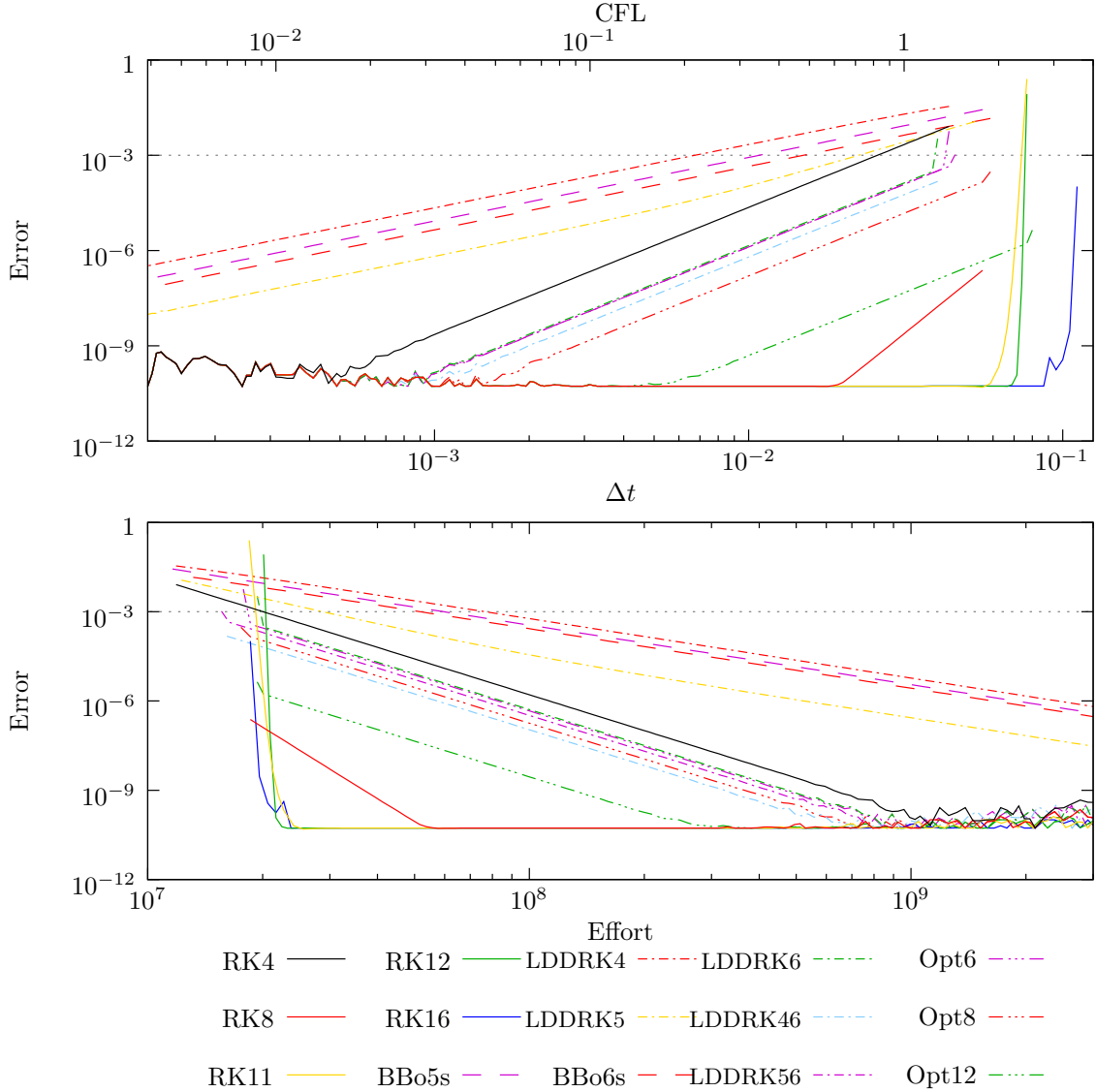}\par%
\caption{Error in the numerical solution of~\eqref{equ:problem} plotted against numerical timestep $\Delta t$ (top plot, bottom scale), or equivalently against $\mathrm{CFL} = \mathrm{PPW}\Delta t$ (top plot, top scale), for various timestepping schemes of varying numerical cost.  All results are using a ``perfect'' 15-point 14th order spatial derivative with $\mathrm{PPW} = 32$ and the ``perfect'' 19-point 16th order $F_{16,4}$ filter, giving an error of around $5\times 10^{-11}$ when used with a ``perfect'' time integration.  The error against computational effort~\eqref{equ:effort} is plotted in the bottom plot.}%
\label{figure:numerics-mo15}%
\end{figure}%
Figure~\ref{figure:numerics-mo15} compares various timestepping schemes for a ``perfect'' 15-point 14th order maximal order spatial derivative with 32 points per wavelength (PPW), using a ``perfect'' spatial filter $F_{16,4}$ at each time step, as described in Ref.~\citenum{brambley-2016-jcp}.  A ``perfect'' time integration would then result in an error of approximately $5\times 10^{-11}$ using this scheme, giving a noise floor due to the spatial discretization used.  As the timestep $\Delta t$, or equivalently the CFL number, is reduced, the error is reduced for each scheme, in general at a rate given by the timestepping scheme's order of accuracy, until this noise floor is reached.  For too large $\Delta t$ the schemes become unstable, generally for CFL numbers in the range 1--4.  The higher order schemes show a significantly lower error than the lower-order schemes, with the 2nd order optimized schemes of~\citet{bogey+bailly-2004} and~\citet[][LDDRK4 and LDDRK5]{hu+hussaini+manthey-1996} performing worse than the 4th order RK4 and optimized LDDRK6, LDDRK46 and LDDRK56~\citep{hu+hussaini+manthey-1996} schemes, which themselves perform worse than the higher order maximal order RK8--16 schemes, although of course the latter involve more stages and therefore a higher computational cost. The Opt12 scheme is able to offer a low error for large $\Delta t$, but this also occurs near the stability limit of the Opt12 scheme. Generally for CFL numbers less than one, the Opt8 and Opt12 offer lower errors than all other optimised scheme, but these schemes also do have significantly more stages than the others. The bottom plot in figure~\ref{figure:numerics-mo15} plots the same error against a measure of the numerical cost of the simulation.  The numerical cost, or numerical effort, is defined to be
\begin{equation}
\textrm{Effort} = pw(T/\Delta t)(L/\Delta x),
\label{equ:effort}  
\end{equation}
where $p$ is the number of Runge--Kutta stages, $w$ is the half-width of the spatial derivative scheme (so the total width is $2w+1$), $T = 24$ is the total simulation time, $L=24$ is the simulation spatial length, and $\Delta t$ and $\Delta x = 1/\mathrm{PPW}$ are the time step and grid spacing.  For a target error of $10^{-3}$, almost all efficient schemes need to be run very close to their stability limit, although this could be an artifact of using an unnaturally accurate spatial discretization.

If we are interested in achieving an overall accuracy of $10^{-3}$, a more conventional spatial discretization would usually be used.
\begin{figure}%
\centering%
\psfrag{Error}[B][B]{Error}%  
\psfrag{dt}[t][t]{$\Delta t$}%  
\psfrag{CFL}[B][B]{CFL}%  
\psfrag{Effort}[t][t]{Effort}%
\psfrag{y1.0e-12}[r][r]{$10^{-12}$}%  
\psfrag{y1.0e-09}[r][r]{$10^{-9}$}%  
\psfrag{y1.0e-06}[r][r]{$10^{-6}$}%  
\psfrag{y1.0e-03}[r][r]{$10^{-3}$}%  
\psfrag{y1.0e-02}[r][r]{$10^{-2}$}%  
\psfrag{y1.0e-01}[r][r]{$10^{-1}$}%  
\psfrag{y2.0e-01}[r][r]{$0.2$}%  
\psfrag{y5.0e-01}[r][r]{$0.5$}%  
\psfrag{y1.0e+00}[r][r]{$1$}%  
\psfrag{xt1.0e-02}[B][B]{$10^{-2}$}%  
\psfrag{xt5.0e-02}[B][B]{$0.05$}%  
\psfrag{xt1.0e-01}[B][B]{$10^{-1}$}%  
\psfrag{xt2.0e-01}[B][B]{$0.2$}%  
\psfrag{xt5.0e-01}[B][B]{$0.5$}%  
\psfrag{xt1.0e+00}[B][B]{$1$}%  
\psfrag{xt2.0e+00}[B][B]{$2$}%  
\psfrag{xt4.0e+00}[B][B]{$4$}%  
\psfrag{x1.0e+09}[B][B]{$10^9$}%  
\psfrag{x1.0e+08}[B][B]{$10^8$}%  
\psfrag{x1.0e+07}[B][B]{$10^7$}%  
\psfrag{x1.0e+06}[B][B]{$10^6$}%  
\psfrag{x5.0e-01}[B][B]{$0.5$}%  
\psfrag{x2.0e-01}[B][B]{$0.2$}%  
\psfrag{x1.0e-01}[B][B]{$10^{-1}$}%  
\psfrag{x5.0e-02}[B][B]{$5\times 10^{-2}$}%  
\psfrag{x2.0e-02}[B][B]{$2\times 10^{-2}$}%  
\psfrag{x1.0e-02}[B][B]{$10^{-2}$}%  
\psfrag{x5.0e-03}[B][B]{$5\times 10^{-3}$}%  
\psfrag{x1.0e-03}[B][B]{$10^{-3}$}%  
\psfrag{RK4}[r][r]{RK4}%  
\psfrag{RK8}[r][r]{RK8}%  
\psfrag{RK11}[r][r]{RK11}%  
\psfrag{RK12}[r][r]{RK12}%  
\psfrag{RK16}[r][r]{RK16}%  
\psfrag{BBo5s}[r][r]{BBo5s}%  
\psfrag{BBo6s}[r][r]{BBo6s}%  
\psfrag{LDDRK4}[r][r]{\small{LDDRK4}}%  
\psfrag{LDDRK5}[r][r]{\small{LDDRK5}}%  
\psfrag{LDDRK6}[r][r]{\small{LDDRK6}}%  
\psfrag{LDDRK46}[r][r]{\small{LDDRK46}}%  
\psfrag{LDDRK56}[r][r]{\small{LDDRK56}}%  
\psfrag{PET6}[r][r]{Opt6}%  
\psfrag{PET8}[r][r]{Opt8}%  
\psfrag{PET12}[r][r]{Opt12}%  
\includegraphics{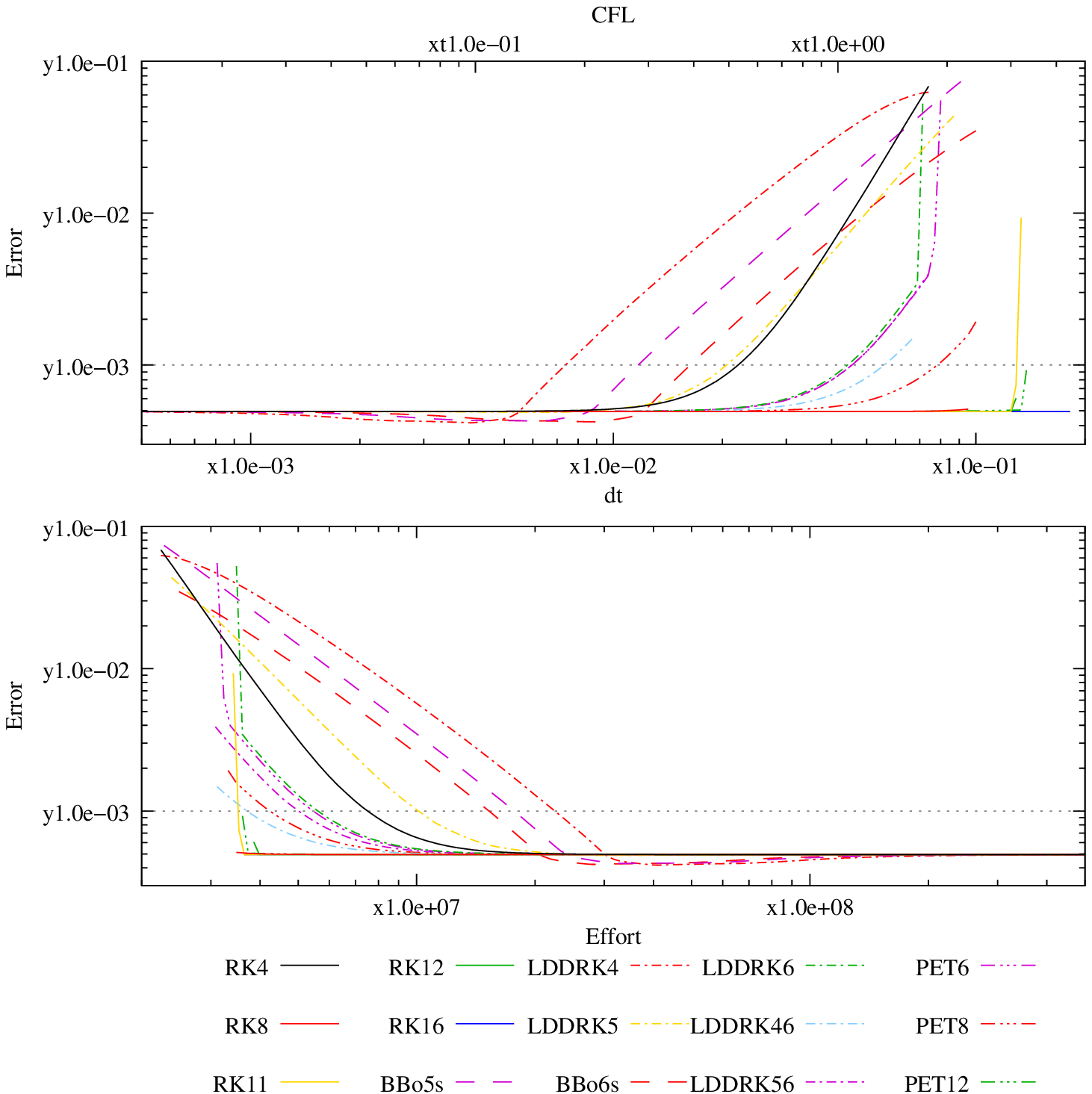}\par%
\caption{Error in the numerical solution of~\eqref{equ:problem} plotted against numerical timestep $\Delta t$ (top plot, bottom scale), or equivalently against CFL (top plot, top scale).  All results are using a 7-point 6th order spatial derivative with $\mathrm{PPW} = 24$ and the 7-point 6th order $F_6$ filter, giving an error of around $5\times 10^{-4}$ when used with a ``perfect'' time integration.  The error against computational effort~\eqref{equ:effort} is plotted in the bottom plot.}%
\label{figure:numerics-mo724}%
\end{figure}
Figure~\ref{figure:numerics-mo724} shows the error for a 7-point 6th~order (maximal order) spatial derivative using 24 PPW and a standard 7-point 6th~order spatial filter.  The noise floor achieved with a ``perfect'' time integration in this case is approximately $5\times 10^{-4}$, which is the limit of accuracy of the spatial discretization.  Apart from the worse performance of the RK4 scheme, the same trend as in figure~\ref{figure:numerics-mo15} is apparent.  In particular, the commonly used LDDRK56 scheme requires $\mathrm{CFL}=1.09$ in order to efficiently achieve the desired accuracy, and the higher order maximal order schemes, while potentially 30\% more computationally efficient for the same accuracy, require an even higher CFL number to achieve this which is close to their stability limit.  The LDDRK46 scheme outperforms the LDDRK56 scheme despite having fewer stages per timestep on average, although needs $\mathrm{CFL}=1.34$ in order to efficiently achieve the desired accuracy; this better performance is possibly due to the larger $\hat{\lambda}_\delta$ for $\delta = 10^{-4}$ seen in table~\ref{table:lddrk} for this scheme, due to its less aggressive optimization. Similarly to the maximal order schemes, the Opt12 scheme rapidly goes to the accuracy limit and any improvements due to the being optimised are not apparent. Also notably, for the same computational error, LDDRK46 outperforms Opt8 which outperforms LDDRK56; this can be attributed to the higher error boundaries for the schemes at $\delta = 10^{-4}$, as seen in table~\ref{tbl:StabilityReoptimisedSchemes}.  However, it could be argued that the spatial discretization in this example has significantly more points per wavelength than is usual in practice for computational aeroacoustics.

\begin{figure}%
\centering%
\psfrag{Error}[B][B]{Error}%  
\psfrag{dt}[t][t]{$\Delta t$}%  
\psfrag{CFL}[B][B]{CFL}%  
\psfrag{Effort}[t][t]{Effort}%
\psfrag{y1.0e-12}[r][r]{$10^{-12}$}%  
\psfrag{y1.0e-09}[r][r]{$10^{-9}$}%  
\psfrag{y1.0e-06}[r][r]{$10^{-6}$}%  
\psfrag{y1.0e-03}[r][r]{$10^{-3}$}%  
\psfrag{y1.0e-02}[r][r]{$10^{-2}$}%  
\psfrag{y1.0e-01}[r][r]{$10^{-1}$}%  
\psfrag{y2.0e-01}[r][r]{$0.2$}%  
\psfrag{y5.0e-01}[r][r]{$0.5$}%  
\psfrag{y1.0e+00}[r][r]{$1$}%  
\psfrag{xt1.0e-02}[B][B]{$10^{-2}$}%  
\psfrag{xt5.0e-02}[B][B]{$0.05$}%  
\psfrag{xt1.0e-01}[B][B]{$10^{-1}$}%  
\psfrag{xt2.0e-01}[B][B]{$0.2$}%  
\psfrag{xt5.0e-01}[B][B]{$0.5$}%  
\psfrag{xt1.0e+00}[B][B]{$1$}%  
\psfrag{xt2.0e+00}[B][B]{$2$}%  
\psfrag{xt4.0e+00}[B][B]{$4$}%  
\psfrag{x1.0e+09}[B][B]{$10^9$}%  
\psfrag{x1.0e+08}[B][B]{$10^8$}%  
\psfrag{x1.0e+07}[B][B]{$10^7$}%  
\psfrag{x1.0e+06}[B][B]{$10^6$}%  
\psfrag{x5.0e-01}[B][B]{$0.5$}%  
\psfrag{x2.0e-01}[B][B]{$0.2$}%  
\psfrag{x1.0e-01}[B][B]{$10^{-1}$}%  
\psfrag{x5.0e-02}[B][B]{$5\times 10^{-2}$}%  
\psfrag{x2.0e-02}[B][B]{$2\times 10^{-2}$}%  
\psfrag{x1.0e-02}[B][B]{$10^{-2}$}%  
\psfrag{x5.0e-03}[B][B]{$5\times 10^{-3}$}%  
\psfrag{x1.0e-03}[B][B]{$10^{-3}$}%  
\psfrag{RK4}[r][r]{RK4}%  
\psfrag{RK8}[r][r]{RK8}%  
\psfrag{RK11}[r][r]{RK11}%  
\psfrag{RK12}[r][r]{RK12}%  
\psfrag{RK16}[r][r]{RK16}%  
\psfrag{BBo5s}[r][r]{BBo5s}%  
\psfrag{BBo6s}[r][r]{BBo6s}%  
\psfrag{LDDRK4}[r][r]{\small{LDDRK4}}%  
\psfrag{LDDRK5}[r][r]{\small{LDDRK5}}%  
\psfrag{LDDRK6}[r][r]{\small{LDDRK6}}%  
\psfrag{LDDRK46}[r][r]{\small{LDDRK46}}%  
\psfrag{LDDRK56}[r][r]{\small{LDDRK56}}%  
\psfrag{PET6}[r][r]{Opt6}%  
\psfrag{PET8}[r][r]{Opt8}%  
\psfrag{PET12}[r][r]{Opt12}%  
\includegraphics{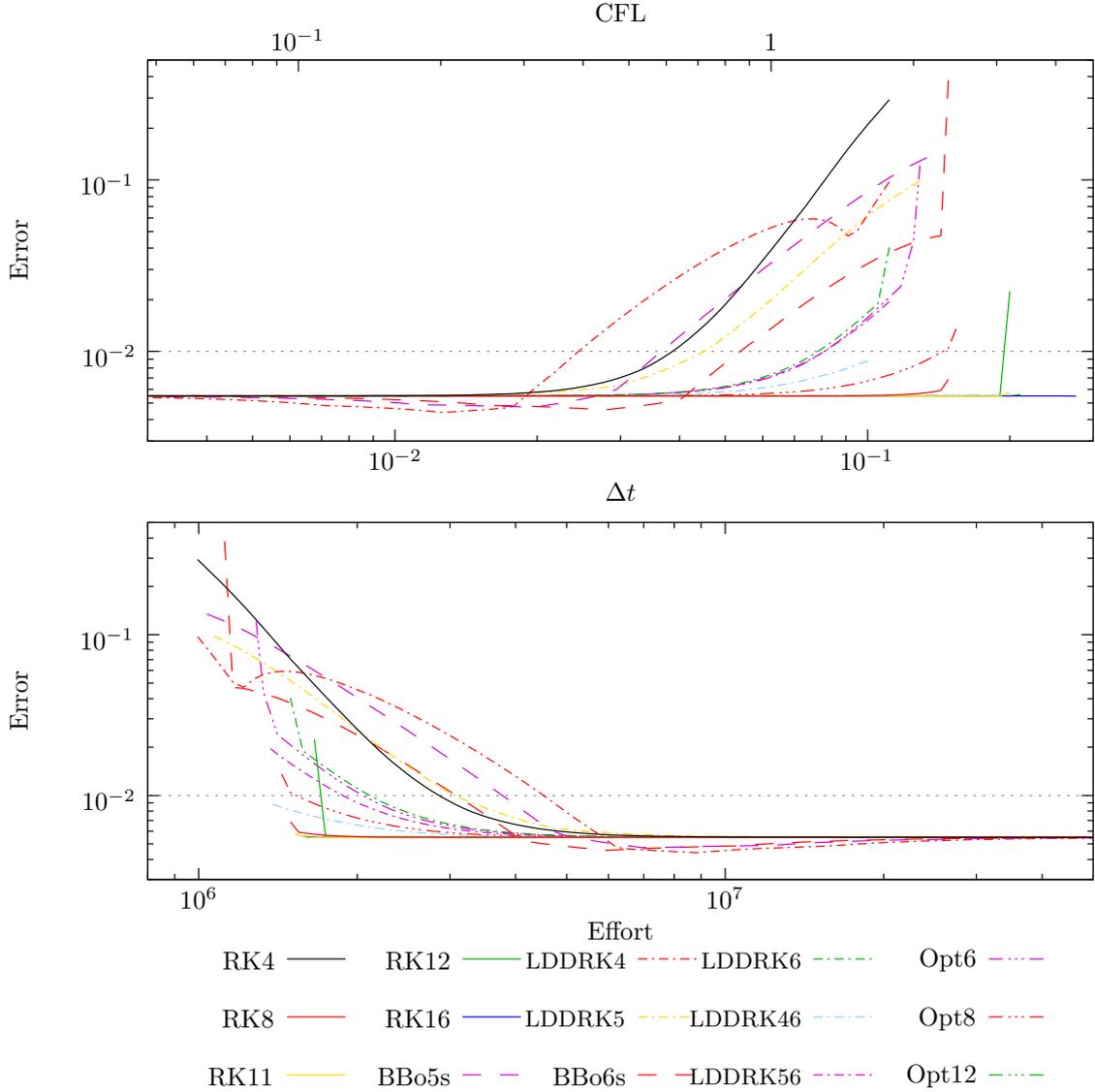}\par%
\caption{Error in the numerical solution of~\eqref{equ:problem} plotted against numerical timestep $\Delta t$ (to plot, bottom scale), or equivalently against CFL (top plot, top scale).  All results are using a 7-point 6th order spatial derivative with $\mathrm{PPW} = 16$ and the 7-point 6th order $F_6$ filter, giving an error of around $5\times 10^{-3}$ when used with a ``perfect'' time integration.  The error against computational effort~\eqref{equ:effort} is plotted in the bottom plot.}%
\label{figure:numerics-mo716}%
\end{figure}%
Dropping the required accuracy to $10^{-2}$ and using 16 points per wavelength results in figure~\ref{figure:numerics-mo716}, with a ``perfect'' time integration noise floor of approximately $5\times 10^{-3}$.  The LDDRK4 scheme starts to show some of the expected optimized behaviour of the optimized timestepping schemes, although at errors of around $0.05$ which do not benefit the desired accuracy of $0.01$.  Once again the higher order maximal order schemes give the best accuracy, and once again they must be run near their stability limit for efficiency.  The most computationally efficient scheme to reach an error of $0.01$ first is the LDDRK46 scheme, although at an unusually high CFL number of $1.6$ at its stability limit. As with figure~\ref{figure:numerics-mo724}, Opt12 almost instantly reaches the noise floor, suggesting it is limited by stability rather than accuracy, similarly to the unoptimized RK8 or RK12 schemes.

In aeroacoustics practice, optimized spatial derivatives are optimized to work efficiently at around 6 points per wavelength or fewer.  In Ref.~\citenum{brambley-2016-jcp} it was shown that such optimized schemes do not perform well in this test case involving non-constant-amplitude waves.
\begin{figure}%
\centering%
\psfrag{Error}[B][B]{Error}%  
\psfrag{dt}[t][t]{$\Delta t$}%  
\psfrag{CFL}[B][B]{CFL}%  
\psfrag{Effort}[t][t]{Effort}%
\psfrag{y1.0e-12}[r][r]{$10^{-12}$}%  
\psfrag{y1.0e-09}[r][r]{$10^{-9}$}%  
\psfrag{y1.0e-06}[r][r]{$10^{-6}$}%  
\psfrag{y1.0e-03}[r][r]{$10^{-3}$}%  
\psfrag{y1.0e-02}[r][r]{$10^{-2}$}%  
\psfrag{y1.0e-01}[r][r]{$0.1$}%  
\psfrag{y2.0e-01}[r][r]{$0.2$}%  
\psfrag{y5.0e-01}[r][r]{$0.5$}%  
\psfrag{y1.0e+00}[r][r]{$1$}%  
\psfrag{xt1.0e-02}[B][B]{$10^{-2}$}%  
\psfrag{xt5.0e-02}[B][B]{$0.05$}%  
\psfrag{xt1.0e-01}[B][B]{$0.1$}%  
\psfrag{xt2.0e-01}[B][B]{$0.2$}%  
\psfrag{xt5.0e-01}[B][B]{$0.5$}%  
\psfrag{xt1.0e+00}[B][B]{$1$}%  
\psfrag{xt2.0e+00}[B][B]{$2$}%  
\psfrag{xt4.0e+00}[B][B]{$4$}%  
\psfrag{x1.0e+09}[B][B]{$10^9$}%  
\psfrag{x1.0e+08}[B][B]{$10^8$}%  
\psfrag{x1.0e+07}[B][B]{$10^7$}%  
\psfrag{x1.0e+06}[B][B]{$10^6$}%  
\psfrag{x5.0e-01}[B][B]{$0.5$}%  
\psfrag{x2.0e-01}[B][B]{$0.2$}%  
\psfrag{x1.0e-01}[B][B]{$0.1$}%  
\psfrag{x5.0e-02}[B][B]{$0.05$}%  
\psfrag{x2.0e-02}[B][B]{$0.02$}%  
\psfrag{x1.0e-02}[B][B]{$0.01$}%  
\psfrag{x5.0e-03}[B][B]{$0.005$}%  
\psfrag{x1.0e-03}[B][B]{$10^{-3}$}%  
\psfrag{RK4}[r][r]{RK4}%  
\psfrag{RK8}[r][r]{RK8}%  
\psfrag{RK11}[r][r]{RK11}%  
\psfrag{RK12}[r][r]{RK12}%  
\psfrag{RK16}[r][r]{RK16}%  
\psfrag{BBo5s}[r][r]{BBo5s}%  
\psfrag{BBo6s}[r][r]{BBo6s}%  
\psfrag{LDDRK4}[r][r]{\small{LDDRK4}}%  
\psfrag{LDDRK5}[r][r]{\small{LDDRK5}}%  
\psfrag{LDDRK6}[r][r]{\small{LDDRK6}}%  
\psfrag{LDDRK46}[r][r]{\small{LDDRK46}}%  
\psfrag{LDDRK56}[r][r]{\small{LDDRK56}}%  
\psfrag{PET6}[r][r]{Opt6}%  
\psfrag{PET8}[r][r]{Opt8}%  
\psfrag{PET12}[r][r]{Opt12}%  
\includegraphics{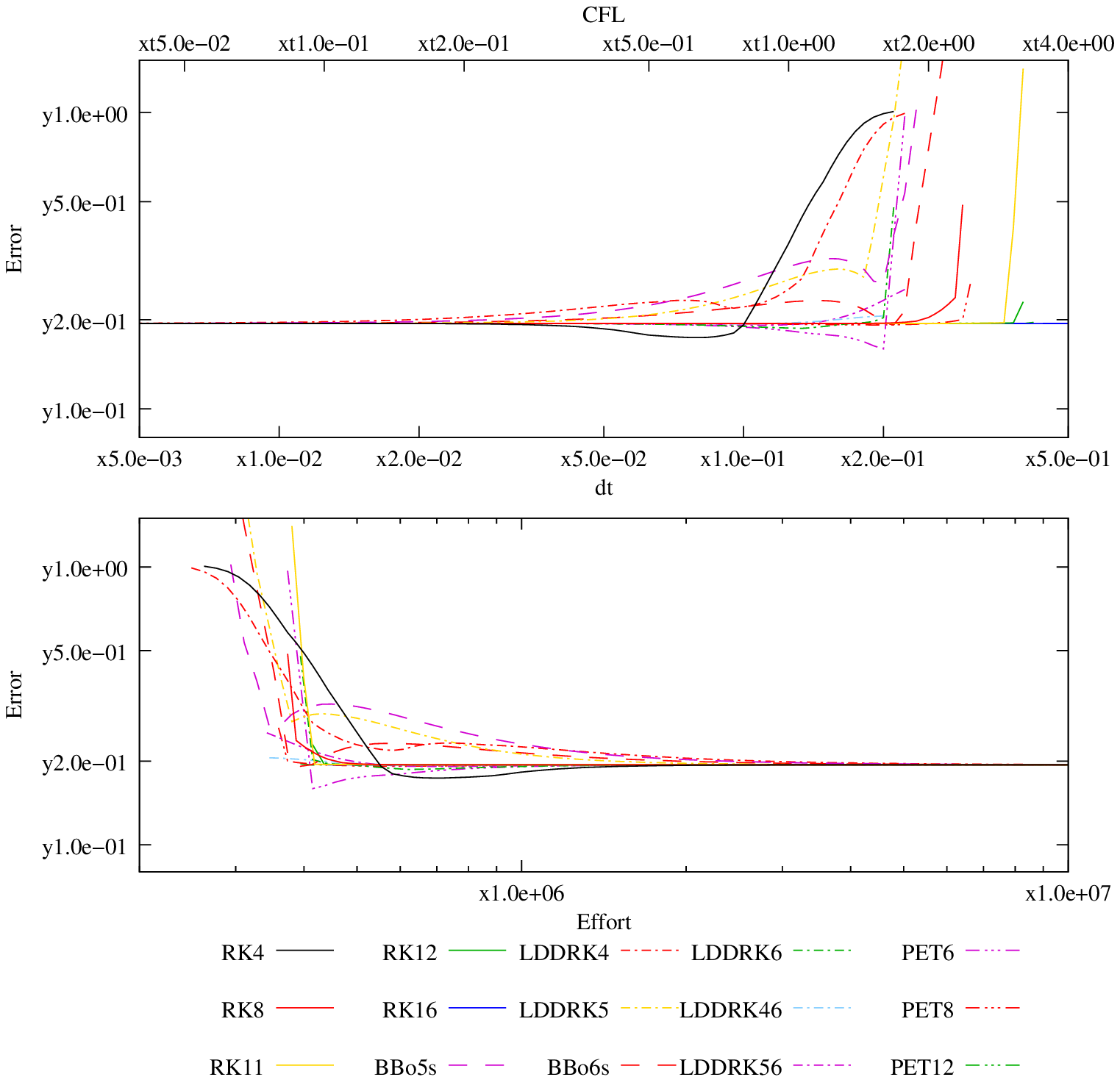}\par%
\caption{Error in the numerical solution of~\eqref{equ:problem} plotted against numerical timestep $\Delta t$ (top plot, bottom scale), or equivalently against CFL (top plot, top scale).  All results are using the 7-point 4th order DRP spatial derivative of~\citet{tam+shen-1993} with $\mathrm{PPW} = 8$ and the 7-point 6th order $F_6$ filter, giving an error of around $0.2$ when used with a ``perfect'' time integration.  The error against computational effort~\eqref{equ:effort} is plotted in the bottom plot.}%
\label{figure:numerics-drp6}%
\end{figure}%
As an example, figure~\ref{figure:numerics-drp6} shows the results of using the 7-points 4th order DRP spatial derivative of~\citet{tam+shen-1993}, together with a standard 7-point 6th order spatial filter.  The noise floor achieved with a ``perfect'' time integration in this case is only $0.2$ (i.e.\ an error of 20\%).  The optimized timestepping schemes can be seen from figure~\ref{figure:numerics-drp6} to also be optimized for this overly ambitious case, with many achieving an error of around $0.3$ or lower more quickly than would be expected from a simple power law error decay. It is unclear what error would be being targeted in this case, and although the optimized schemes do seem to outperform the maximal order schemes when comparing error against computational effort in this case, the errors are all sufficiently large that no scheme could be said to have properly converged.

%%%%%%%%%%%%
%Conclusion%
%%%%%%%%%%%%

\section{Conclusion}
\label{section:conclusion}

This paper has investigated Runge--Kutta timestepping schemes optimized to solve linear time-invariant oscillatory problems.  By Fourier transforming, without loss of generality we have focused attention on oscillations with time dependence $\exp\{-\I\omega t\}$.  By analogy with optimized spatial derivatives (DRP schemes), which were found to behave poorly for waves of growing or decaying amplitudes~\citep{brambley-2016-jcp}, it is found here that existing optimized timestepping schemes have also been optimized assuming real frequencies $\omega$, and hence assuming constant amplitude oscillations, and that these schemes also perform poorly for oscillations of growing or decaying amplitude, corresponding to complex values of $\omega$.  In particular, figure~\ref{fig:twoStepLDDRKvstwoTimeStepRkOnePercentPhaseError} shows that maximal order schemes are more accurate than optimized LDDRK46 or LDDRK56 schemes~\citep{hu+hussaini+manthey-1996} for the same computational effort apart from for a rather limited range of nearly-real frequencies; targeting these frequencies would require significant a priori knowledge of the simulation to be performed, and is likely impossible for broadband simulations.

Section~\ref{section:optimization} attempts to reoptimize Runge--Kutta schemes by considering their behaviour throughout the complex $\omega$ plane, and not just for real $\omega$.  A wide range of different optimizations were attempted, and a selection of three well-performing schemes is given in table~\ref{table:OptimizedCoefficients}.  These three schemes, labelled here Opt6, Opt8, and Opt12, are reasonably typical of the types of schemes obtained from such optimizations, and, while our search for such schemes was not exhaustive, we believe were there significantly more accurate schemes they would have been found.  Unsurprisingly, these new schemes perform better for complex values of $\omega$ than the existing schemes optimized over only real $\omega$.  However, the newly optimized schemes only outperform the more traditional maximal order schemes for a relatively small region of complex $\omega$, and as commented above, targeting these frequencies would require significant a priori knowledge of the simulation to be performed.  From theoretical considerations, therefore, we conclude that either an 8- or 12-stage maximal order Runge--Kutta scheme would appear to perform best for broadband simulations of linear non-constant-amplitude oscillations, with the timestep dictated by stability rather than accuracy considerations.

The theory was illustrated by solving an example 1D wave equation in section~\ref{section:1d} taken from Ref.~\citenum{brambley-2016-jcp}.  The results suggest that optimized timestepping schemes have been over-ambitiously optimized with target errors of around $\delta = 10^{-3}$ \emph{per timestep}, where as significantly smaller errors are required to get an overall simulation error of around $10^{-2}$ or $10^{-3}$.  The results also suggest that the best computational efficiency for a given target accuracy is obtained when the timestepping scheme is very close to its stability limit, with unusually large CFL numbers of 1.5--2 being typically optimal.  As was found in Ref.~\citenum{brambley-2016-jcp}, significantly more points per wavelength are needed for the spatial derivatives (around 16 PPW) than are commonly thought necessary to achieve even a modest accuracy of 1\% error, and optimized DRP spatial derivatives perform worse than maximal order spatial derivatives.  It is notable that the LDDRK46 scheme~\citep{hu+hussaini+manthey-1996} outperforms the more widely used LDDRK56 scheme in all cases considered here.  However, for well-resolved oscillations, the example 1D wave equation gives the most accurate results with the maximal order RK8 or RK12 schemes, with the error being limited by the spatial derivatives in all cases.

While the accuracy of timestepping schemes for non-constant-amplitude waves (with complex $\omega$) is a straightforward extension of the notion of accuracy for constant-amplitude waves (with real $\omega$), the same is not true for the stability of timestepping schemes.  Indeed, a timestepping scheme is stable for constant-amplitude waves if the numerical solution does not grow in time, meaning only under-predicted growth rates are allowed.  It is probably undesirable to require that the amplitude of non-constant-amplitude waves is always under-predicted numerically, and instead growth- and decay-rates of non-constant-amplitude waves are desired to be modelled numerically as accurately as possible.  This suggests that the concept of stability for timestepping schemes is restricted to only constant-amplitude waves (with real $\omega$), and no extension to the concept of stability is needed for non-constant-amplitude waves.

It is worth noting that this work has solely considered \emph{linear} oscillatory problems, for which a maximal order $p$-stage Runge--Kutta timestepping scheme can achieve $p$th order accuracy.  All schemes considered here will only be 2nd order accurate when applied to nonlinear problems.  In such cases, the standard RK4 scheme retains 4th order accuracy even for nonlinear problems.  To our knowledge, the optimization of such higher order schemes for nonlinear problems has not been attempted.  We mention in passing that a linear wave equation with a nonzero source term classifies as a nonlinear problem in this context.

Another possibility for future investigation would be the joint optimization of spatial derivatives and temporal time-stepping schemes, as was performed under the assumption of constant-amplitude oscillations by \citet{rona+etal-2017}.  However, such optimizations would need to assume an underlying dispersion relation linking the spatial wavenumbers and temporal frequencies, and such optimizations are likely to be specific to the particular dispersion relation chosen.  This is unlike the optimization of just temporal or spatial derivatives in isolation, which should be generally applicable to systems with arbitrary dispersion relations.

\section*{Acknowledgements}

AP gratefully acknowledges funding from the Undergraduate Research Summer Studentship (URSS) scheme at the University of Warwick.
EJB gratefully acknowledges funding from a Royal Society University Research Fellowship (UF150695).
\prevpub

%%%%%%%%%%%%
%%Appendix%%
%%%%%%%%%%%%

\appendix
\section{Stability Derivation}\label{app:Stability}

In this appendix we derive the conditions for a Runge--Kutta scheme to be stable in the limit $\odt\to0$.  For any timestepping scheme, solving the differential equation $\intd\vect{U}/\intd t = -\I\omega\vect{U}(t)$ gives the solution $\vect{U}(t+\Delta t) = \vect{U}(t)r(\odt)$ where $r(\odt)$ is the amplification factor.  Stability is given by the amplification factor having a modulus less than one.  Any $q$th order $p$-stage explicit Runge--Kutta scheme has an amplification factor of the form
\begin{align}
r(\odt) &= 1 + \sum_{j=1}^q \frac{1}{j!}(-\I\odt)^j + \sum_{\mathclap{j=q+1}}^p c_j(-\I\odt)^j,
\end{align}
for some coefficients $c_j$.  The sum up to $p$ may be extended to a sum to infinity by defining $c_j = 0$ for $j>p$.  The amplification factor for a perfectly accurate scheme would be $r_e(\odt) = \exp\{-\I\odt\}$.  Noting that $|r_e(\odt)|=1$, we may write the stability condition as $|r(\odt)|<1 \Leftrightarrow \Real(\log(r(\odt)/r_e(\odt))) < 0$.  Expanding for small $\odt$,
\begin{subequations}\begin{align}
\log\left(\frac{r(\odt)}{r_e(\odt)}\right) =& \sum_{j=1}^q \frac{1}{j!}(-\I\odt)^j + \sum_{\mathclap{j=q+1}}^\infty c_j(-\I\odt)^j - \frac{1}{2}\!\left(\sum_{j=1}^q \frac{1}{j!}(-\I\odt)^j + \sum_{\mathclap{j=q+1}}^\infty c_j(-\I\odt)^j\right)^{\!\!2}
\\\notag&
- \left(\sum_{j=1}^q \frac{1}{j!}(-\I\odt)^j + \!\!\!\sum_{j=q+1}^\infty \frac{1}{j!}(-\I\odt)^j\right)\! + \frac{1}{2}\!\left(\sum_{j=1}^q \frac{1}{j!}(-\I\odt)^j + \!\!\!\sum_{j=q+1}^\infty \frac{1}{j!}(-\I\odt)^j\right)^{\!\!2}
\\\notag&+ \cdots\\
=&\left(c_{q+1}-\frac{1}{(q+1)!}\right)\!(-\I\odt)^{q+1}
\\\notag&\qquad\qquad
+ \left[\left(c_{q+2}-\frac{1}{(q+2)!}\right) - \left(c_{q+1}-\frac{1}{(q+1)!}\right)\right]\!(-\I\odt)^{q+2}
+ O\big((\odt)^{q+3}\big).
\end{align}\end{subequations}
Hence, taking the real part, for $\omega$ being real,
\begin{align}
&\Real\!\left(\log\left(\frac{r(\odt)}{r_e(\odt)}\right)\!\right) =
\\\notag&\qquad
\left\{\begin{array}{ll}
  (-1)^{n+1}\left[\left(c_{2n+2}-\frac{1}{(2n+2)!}\right) - \left(c_{2n+1}-\frac{1}{(2n+1)!}\right)\right]\!(\odt)^{2n+2} + O\big((\odt)^{2n+4}\big)
  & \text{if $q=2n$,}\\
  (-1)^n\left(c_{2n}-\frac{1}{(2n)!}\right)\!(\odt)^{2n} + O\big((\odt)^{2n+2}\big)
  & \text{if $q=2n-1$,}
  \end{array}\right.
\end{align}
This is the result given in~\eqref{equ:instability}.  Requiring this to be less than zero gives the stability conditions following~\eqref{equ:instability}.

\ifAIAA\else\section*{Bibliography}\fi

\bibliography{paper}

\end{document}